\RequirePackage{fix-cm}
\documentclass[smallextended]{svjour3}       
\smartqed  
\usepackage{graphicx}
\usepackage[title]{appendix}
\usepackage{amssymb,color,ulem,amsmath,cancel}

\newcommand{\sign}{{\rm \hskip0.5pt sign \hskip1pt}}

\def\ra{\rangle}
\def\la{\langle}

\usepackage{cases}
\usepackage{soul}

\begin{document}

\title{Stabilization of the response of cyclically loaded lattice spring models with plasticity 
}

\titlerunning{Stabilization of elastoplastic systems}        

\author{Ivan Gudoshnikov         \and
        Oleg Makarenkov
}

\institute{I. Gudoshnikov\at
              University of Texas at Dallas, 
800 West Campbell Road, 
Richardson, TX 75080  \\
              \email{ixg140430@utdallas.edu}           
           \and
           O. Makarenkov \at
              University of Texas at Dallas, 
800 West Campbell Road, 
Richardson, TX 75080 \\
              Tel.: +1-972-883-4617\\
              \email{makarenkov@utdallas.edu}           
}

\date{Received: date / Accepted: date}

\maketitle

\begin{abstract} 
	This paper develops an analytic framework to design both stress-controlled and displacement-controlled $T$-periodic loadings which make the quasistatic evolution of a one-dimensional network of elastoplastic springs converging to a unique periodic regime. The solution of such an evolution problem is a function $t\mapsto (e(t),p(t))$, where $e_i(t)$ and $p_i(t)$ are the elastic and plastic deformations of spring $i,$  defined on $[t_0,\infty)$ by the initial condition $(e(t_0),p(t_0))$. 
 After we rigorously convert the problem into a Moreau sweeping process with a moving polyhedron $C(t)$ in a vector space $E$ of dimension $d,$ it becomes natural to expect (based on a result by Krejci) that the solution $t\mapsto (e(t),p(t))$ always converges to a $T$-periodic function. The achievement of this paper is in spotting a class of loadings where the Krejci's limit doesn't depend on the initial condition $(e(t_0),p(t_0))$ and so all the trajectories approach the same $T$-periodic regime. The proposed class of sweeping processes is the one for which the normals of any $d$ different facets of the moving polyhedron $C(t)$ are linearly independent. 
We further link this geometric condition to mechanical properties of the given network of springs. We discover that
 the normals of any $d$ different facets of the moving polyhedron $C(t)$ are linearly independent, if the number of displacement-controlled loadings is two less the number of nodes of the given network of springs and when the magnitude of the stress-controlled loading is sufficiently large (but admissible). The result can be viewed as an analogue of the high-gain control method for elastoplastic systems. In continuum theory of plasticity, the respective result is known as Frederick-Armstrong theorem. 

\vskip0.2cm

\noindent The theoretical results are accompanied by analytic computations for instructive examples. In particular, we convert a specific one-dimensional network of elastoplastic springs into sweeping process which have never been explicitly addressed in the literature so far.

\end{abstract}

\setcounter{secnumdepth}{2} \setcounter{tocdepth}{2}

\sloppy

\section{Introduction}

The classical theory of elastoplasticity  offers comprehensive results, commonly known as shakedown theorems,  about the maximal magnitude of the applied loading ({\it shakedown load limit}) beyond which the response of elastoplastic material necessarily involves plastic deformation regardless of the initial distribution of stresses in the material, see  \cite[\S10]{jirasek}. In other words, shakedown theorems measure the distance between the current stress distribution in the material to a certain boundary (called {\it yield surface}) built of the spatially distributed elastic limits. The fundamental result by Frederick and Armstrong \cite{frederick} says that, if the amplitude of  a $T$-periodic loading exceeds the shakedown limit, then the stress distribution asymptotically approaches a unique $T$-periodic steady cycle which doesn't depend on the initial stress distribution ({\it uniqueness of the response}).  Frederick-Armstrong highlight that convergence to a unique cyclic state is guaranteed when the yield surface contains no lines of zero curvature \cite[p.~159]{frederick}. Assuming this or another restriction on the geometry of the yielding surface  (such as von Mises, Tresca, or Mohr-Coulomb criteria), many authors  computed the steady cycle by discretizing the problem spatially \cite{garcea,Heitzer,li} and/or temporarily \cite{polizzotto,ponter,brazil},  and by solving the associated minimization problems for the successive discrete states. Applications included the performance of various structures and metal matrix composites under  cyclic loadings, see  \cite{simplified,weichert1}.

\vskip0.2cm

\noindent Aiming to design materials with better properties, there has been a great deal of work lately where a discrete structure comes not from an associated model of continuum mechanics, but from a certain microstructure formulated through a lattice of elastic springs 
\cite{crystal1,li1} (metals),  \cite{polymer,novel} (polymers), \cite{JiaoPrior2} (titanium alloys),  \cite{tissue,natureb,sva} (biological materials). 
Despite of the fact that  fatigue crack initialization in heterogeneous materials strongly depends on local micro-plasticity (see e.g. Blechman \cite{blechman}), the current literature features only numeric results about the dynamics of the lattices of elastoplastic springs. Important papers in this direction are e.g. Buxton et al \cite{buxton} and Chen et al \cite{chen}. 

\vskip0.2cm

\noindent The goal of the present paper is to initiate the development of a qualitative theory of the lattices of elastoplastic springs and to offer an analogue of the Frederick-Armstrong theorem for such systems.

\vskip0.2cm

\noindent We stick to the setting of ideal plasticity (the stress of each spring is constrained within so-called {\it elastic limits} beyond which plastic deformation begins) and investigate the asymptotic distribution of the stresses $s(t)=(s_1(t),...,s_m(t))$ of a network of $m$ elastoplastic springs. Starting with a graph of $m$ connected elastoplastic springs, the paper takes Moreau's approach \cite{moreau} to write down the equation for stress $s_i$ of spring $i$ without any knowledge about plastic deformation of the spring and relying entirely on the geometry of the graph and elastic limits $[c_i^-,c_i^+]$ of the springs. The plasticity is accounted through the $m$-dimensional parallelepiped-shaped constraint $C(t)$,
whose boundary can be viewed as {\it discretized yield surface}, see Fig.~\ref{fff}.
 Beneficially for the performance of computational routines, Moreau concluded that the stress-vector $s(t)=(s_1(t),...,s_m(t))$ of springs is confined within a time-independent low-dimensional hyperplane $V$. 
 It is also due to Moreau that external time-varying loadings enter the equations of dynamics through a time-varying vector $c(t)$ that acts as displacement of the parallelepiped $C$ (Fig.~\ref{fff}). The only obstacle towards practical implementation of the Moreau approach \cite{moreau} in the context of  spring network modeling is that \cite{moreau} deals with abstract configuration spaces translated into practical quantities only for examples from dry friction mechanics. This paper clears this obstacle and fully adapts Moreau sweeping process framework to the modeling of networks of elastoplastic springs. 
 
 \vskip0.2cm
 
\noindent After a suitable change of variables that we rigorously incorporate in the next section of the paper, the equations of Moreau ({\it Moreau sweeping process}) can be formulated as 
\begin{equation}
-y'(t)\in N_{({C}+c(t))\cap V}(y(t)),\quad y(t)\in\mathbb{R}^m,
\label{eq:MSP0}
\end{equation}
where 
\begin{equation}\label{NC}
  N_{\mathcal{C}}(x)=\left\{\begin{array}{ll}\left\{\zeta \in \mathbb{R}^n:\langle\zeta,c-x\rangle\leqslant 0,\ {\rm for\ any }\ c\in {\mathcal{C}}\right\},& {\rm if}\ x\in {\mathcal{C}},\\
   \emptyset,& {\rm if}\ x\not\in {\mathcal{C}}
\end{array}\right.
\end{equation}
is a {\it normal cone} to the set $\mathcal{C}$ at point $x$ and $V$ is a subspace of $\mathbb{R}^m$. For Lipschitz-continuous $t\mapsto c(t)$ (which we show to be the case when the external loading is Lipschitz-continuous) sweeping process (\ref{eq:MSP0}) possesses usual properties of the existence and continuous dependence of solutions on the initial conditions, see e.g. Kunze and Monteiro Marques \cite{kunze}.


 \vskip0.2cm

\noindent When sweeping process (\ref{eq:MSP0}) includes a vector field on top of the normal cone (so-called {\it perturbed sweeping process}), multiple results are available to stabilize the dynamics of a sweeping process. Important results in this direction are obtained in
Leine and van de Wouw \cite{leine,leine2}, Brogliato \cite{brogliato1}, and Brogliato-Heemels  \cite{brogliato2}, Kamenskiy et al \cite{nahs}. 

\vskip0.2cm

 \noindent  As for the regular sweeping process (\ref{eq:MSP0}), very limited tools to control the asymptotic response are currently available (in contrast to optimal control results developed e.g. in Colombo et al \cite{colombo}).  The asymptotic behavior of sweeping process (\ref{eq:MSP0}) with $T$-periodic excitation $t\mapsto c(t)$ was studied in Krejci \cite{Krejci1996}, who proved the convergence of solutions of (\ref{eq:MSP0}) to a $T$-periodic attractor in the case $V=\mathbb{R}^m$, i.e. $(C+c(t))\cap V=C+c(t).$ If sweeping process (\ref{eq:MSP0}) decomposes into a cross-product of several sweeping processes (\ref{eq:MSP0}) with $m=1$, the global asymptotic stability can be concluded from the theory 
of Prandtl-Ishlinskii operators (Brokate-Sprekels \cite{brokate}, 
Krasnosel'skii-Pokrovskii \cite{pokrov}, Visintin \cite{visintin}).
In the case of an arbitrary $T$-periodic polyhedron $t\mapsto(C+c(t))\cap V$, 
it looks possible to follow the ideas of Adli et al \cite{adly} and obtain global asymptotic stability of a periodic solution by assuming that $g(t)$ lies strictly inside the normal cone $N_{(C+c(t))\cap V}(y)$ for at least one $y$ and $t\in[0,T]$. The present paper takes a different route  
and establishes convergence of  solutions to a unique $T$-periodic regime in terms of the shape of the moving constraint $(C+c(t))\cap V$ only.

 \vskip0.2cm

 \noindent The paper is organized as follows. The next section rigorously formulates the system of laws of quasistatic evolution for a one-dimensional network of $m$ elastoplastic springs on $n$ nodes.
 In section~3 we construct  the vector $c(t)$ and the hyperplane $V$ for arbitrary networks of elastoplastic springs of 1-dimensional nodes. We  discover that the functions $g(t)\in V$ and $h(t)\perp V$ in the orthogonal decomposition of $c(t)$ (see Fig.~\ref{fff}) correspond to {\it displacement-controlled}  and {\it stress-controlled} loadings respectively (as termed in \cite{simplified}). 
The achievement of Section~2 makes it possible to link  the dynamics of networks of elastoplastic springs to the dynamics of sweeping processes.


\vskip0.2cm

\noindent In Section~4 we consider a general sweeping process with a moving set of a form $\cap_{i=1}^k(C_i+c_i(t))$, where $C_i$ are closed convex sets, and prove (Theorem~\ref{thm1}) the convergence of all solutions to a $T$-periodic attractor $X(t)$. Section~\ref{strengthening} (Theorem~\ref{thmnew}) sharpens the conclusion of  Theorem~\ref{thm1} for the case when $\cap_{i=1}^k(C_i+c_i(t))$ is the polyhedron $\Pi(t)\cap V.$ Theorem~\ref{thmnew} shows that even though $X(t)$ may consist of a family of functions, all those functions exhibit certain similar dynamics. Specifically, we prove that any two function $x_1,x_2\in X$ reach (leave) any of the facets of   
$\Pi(t)\cap V$ at the same time.
 Section~\ref{sec1a} (Theorem~\ref{cor0}) reformulates the conclusion of Theorem~\ref{thmnew} in terms of the sweeping process of a network of elastoplastic springs.  

\vskip0.2cm

\noindent Section~\ref{sec5} introduces a class of networks of elastoplastic springs whose stresses converge to a unique $T$-periodic regime regardless of applied $T$-periodic loadings as long as the magnitudes of those loadings are sufficiently large. We begin Section~\ref{sec5} by addressing a general sweeping process in a vector space $E$ of dimension $d$ with a $T$-periodic polyhedral moving set with no connection to networks of springs. Theorem~\ref{thm1.5} of Section~\ref{sec51} states  that the periodic attractor of such a sweeping process contains at most one non-constant solution, if 
normals of any $d$ different facets of the moving polyhedron $C(t)$ are linearly independent. 
 Section~\ref{sec52}  is the main achievement of this paper, where we introduce a class of networks of elastoplastic springs for which the condition of Theorem~\ref{thm1.5} can be easily expressed in terms of the magnitudes of the periodic loadings. We discovered (Theorem~\ref{mainthm}) that global stability of a unique periodic regime occurs when both displacement-controlled and stress-controlled loadings are large enough.



\section{The laws of quasistatic evolution for one-dimensional networks of elastoplastic springs}

We consider a one-dimensional network of $m$ elastoplastic springs with elongations $e_k+p_k,$ $k\in \overline{1,m}$, where  $e_k$ and $p_k$ are elastic and plastic components respectively. 
The bounds of the stress of spring $k$ are denoted by $[c_k^-,c_k^+]$ and $a_k$ stays for the Hooke's coefficient of this spring. 
Each spring
 connects two of $n$ nodes according to  $e_k+p_k=
\xi_{j_k}-\xi_{i_k}$, where $i_k$ and $j_k$ are the indices of the left and right nodes of spring $k$ respectively and $\xi_i$ is the displacement of node $i.$ So defined, the one-dimensional network of springs is an oriented graph on $n$ nodes, where the direction from $i_k$ to $j_k$ is viewed positive through $k\in\overline{1,m}.$

\vskip0.2cm

\noindent The paper investigates the evolution of the stresses under the influence of two types of loadings being displacement-controlled loading and stress-controlled loading.

\subsection{Displacement-controlled loading}
\label{ssect:DPL}
\begin{figure}[h]\center
\vskip-0.5cm
\includegraphics{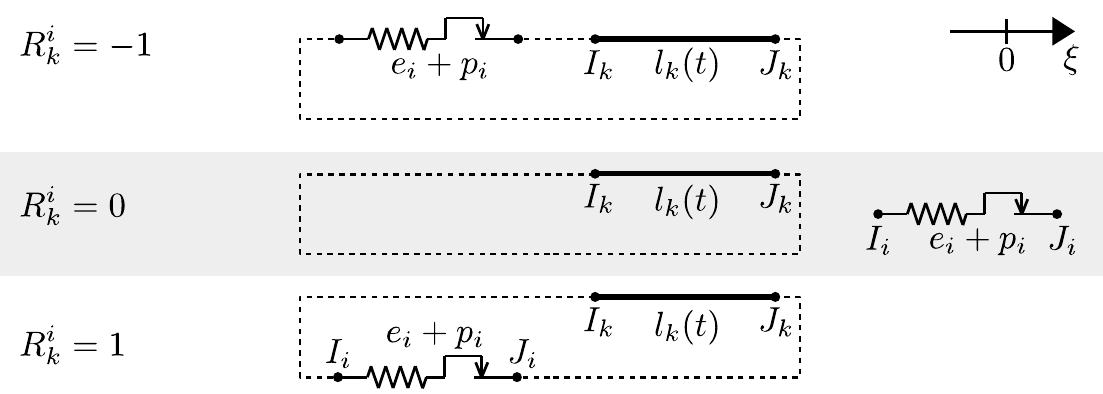}
\caption{\footnotesize Illustration of the signs of the components of the incidence vector $R^k\in\mathbb{R}^m.$ The dotted contour stays for the chain of the springs associated with the vector $R^k$. 
} \label{figR}
\end{figure}

\noindent Displacement-controlled loading  locks the distance between nodes $I_k$ and $J_k$ through $k\in\overline{1,q}$ according to $\xi_{J_k}-\xi_{I_k}=l_k(t)$. Since we will work with {\it connected graphs} of springs only, we assume that 
each length $l_k$ is uniquely determined by the lengths of springs, i.e. 
for each displacement-controlled loading $k\in\overline{1,q}$ there exists a chain of springs which connects the left node $I_k$ of the constraint $k$ with its right node $J_k.$
To each displacement-controlled loading  $k$ we can, therefore, associate a so-called {\it incidence vector}  $R^k\in\mathbb{R}^m$ whose $i$-th component $R^k_i$ is $-1,$ $0,$ or $1$ according to whether the spring $i$ increases, not influences, or decreases the displacement when moving from node $I_k$ to $J_k$ along the chain selected, see Fig.~\ref{figR}


\subsection{Stress-controlled loading} \noindent The stress-controlled loading $f_i$ models an external force applied at node $i,$ so that it affects the resultant of forces at node $i.$  We will study a so-called {\it quasistatic evolution} problem which further assumes that 
 $f(t)$ can be balanced by the stresses of springs at any time. In other words, we assume that the stresses of springs, the reactions of displacement-controlled constraints, and the applied stress loading compensate one another at each of the $n$ nodes. In particular, for any $f(t)=const$ the system admits an  equilibrium.

\vskip0.3cm

  \subsection{The variational system}
\begin{figure}[h]\center
\vskip-0.2cm
\includegraphics{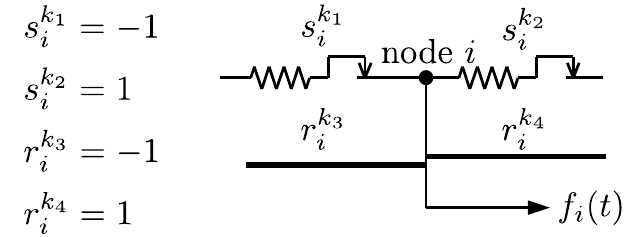}\vskip-0.0cm
\caption{\footnotesize Examples of forces applied at node $i$.
} \label{staticlaw}
\end{figure}

\noindent With the notations introduced the quasistatic evolution of the stresses $s_k$ of springs and reactions  $r_k$ of displacement-controlled loadings can be described by the following variational system (which corresponds to equations (6.1)-(6.6) in the abstract framework by Moreau \cite{moreau})
\begin{eqnarray}
\hskip-0.5cm \mbox{Elastic deformation:} && s=Ae,\label{eq1}\\
\mbox{Plastic deformation:} 
&&  \dot p\in N_C(s), \label{formula2}\\
\mbox{Geometric constraint:}&&   e+p\in D\mathbb{R}^n,\label{D}\\
\mbox{Displacement-controlled loading:} &&R^\top(e+p)=l(t),\label{G}\\
   \mbox{Static balance under\ \ \ \ \ \ }  && s_ks^1+\ldots+s_ms^m+\nonumber\\
\mbox{stress-controlled loading:}&& +r_1r^1+\ldots+r_qr^q+f(t)=0,\hskip1cm \label{eq5}
\end{eqnarray}
where
\begin{eqnarray*}
  s=(s_1,\ldots s_m)^\top &-& \mbox{stresses of springs,}\\
  r=(r_1,\ldots r_q)^\top &-& \mbox{reactions of displacement-controlled loadings,}\\
 e=(e_1,\ldots e_m)^\top &-& \mbox{elastic elongations of springs},\\ p=(p_1,\ldots p_m)^\top &-& \mbox{plastic elongations of springs,}\\
  l(t)=\left(l_1(t),\ldots,l_q(t)\right)^\top &-& \mbox{enforced \hskip-0.04cm lengths \hskip-0.04cm  between \hskip-0.04cm  nodes }\hskip-0.04cm  I_i\mbox{ \hskip-0.04cm  and   } \hskip-0.04cm J_i,\hskip-0.04cm  \ i\in\overline{1,q},\\
f(t)=\left(f_1(t),\ldots,f_n(t)\right)^\top &-& \mbox{stress-controlled loadings at nodes,}\\
A={\rm diag}(a_1,...,a_m) &-& \mbox{matrix of Hooke's coefficients,} \\
N_C(s)=\otimes_{i=1}^m N_{[c_i^-,c_i^+]}(s) &-& \text{the normal cone to } C=\otimes_{i=1}^m [c_i^-,c_i^+]\mbox{ at }s,  \\
D\xi=\left(\xi_{j_k}-\xi_{i_k}\right)_{k=1}^m &-& \text{a linear map (represented by } m\times n\text{-matrix $D$)}\\&& \text{that defines the graph of springs} ,\\ 
R=\left(R^1,\ldots,R^q\right)&-&  m\times q\mbox{-matrix  of the incidence vectors},\\ 
& &  \mbox{of displacement-controlled loadings}.
\end{eqnarray*}
while  the vectors $s^k=(s^k_1,\ldots s^k_n)^\top$ and $r^k=(r^k_1,\ldots r^k_n)^\top$ describe the signs of contributions of stresses of spring $k$ and reactions of displacement-controlled loading  $k$ into the resultant of forces at nodes $1,...,n,$ i.e. (see Fig.~\ref{staticlaw})
\begin{itemize}
\item
$s^k_i=-1,\ s^k_i=0,\ {\rm or}\ s^k_i=1,$  according to whether the spring $k$ is to the left from node $i,$  not connected to node $i,$  or is to the right from node $i,$
\item
$r^k_i=-1,\ r^k_i=0,\ {\rm or}\ r^k_i=1,$ according to whether the displacement-controlled loading $k$  is applied to the left from node $i,$  not applied to node $i,$ or applied to the right from node $i.$
\end{itemize}
\noindent The $m\times n$-matrix $D$ will be termed the {\it kinematic matrix} of  the one-dimensional network of $m$ springs on $n$ nodes. Note, the matrix $-D^\top$ will then be the incidence matrix of the associated oriented graph of $n$ nodes and $m$ edges. 
\vskip0.2cm

\noindent The static balance law \eqref{eq5} will further be written in the equivalent shorter form (\ref{Ivan11}), which is similar to the one used by Moreau \cite[formula (3.23)]{moreau}. We believe, however, that formulation (\ref{eq5}) creates a better idea as for why this law does indeed balance the forces.

\vskip0.2cm

\noindent Following Moreau \cite{moreau}, we term system (\ref{eq1})-(\ref{eq5}) an {\it elastoplastic system.}

\section{Casting the variatonal system as a sweeping process}\label{sec3}
\subsection{Derivation of the sweeping process}
\noindent In order for (\ref{G}) to be solvable in $e+p$  we assume that the displacement-controlled loadings $\{l_i(t)\}_{i=1}^q$ are independent in the sense that  
\begin{equation}\label{A2}
{\rm rank}\left(D^\top R\right)={\rm rank}\left(R^\top D\right)=q.
\end{equation}
Mechanically, condition (\ref{A2}) ensures that the displacement-controlled loadings don't contradict one another. For example, (\ref{A2}) rules out the situation where two displacement-controlled loadings connect same pair of nodes. 
It follows from condition \eqref{A2} that the  matrix equation 
\begin{equation}
   R^\top D L=I_{q\times q}.
   \label{barxi}
\end{equation}
has a $n\times q$ matrix solution. Furthermore, as we will show in the proof of Theorem~\ref{moreauthm}, in order for equation (\ref{eq5}) to be solvable in $s\in\mathbb{R}^m$ and $r\in\mathbb{R}^q$, the function $f(t)$ must satisfy $f(t)\in D^\top\mathbb{R}^m$. That is why, the existence of a continuous function $\bar h:\mathbb{R}\to\mathbb{R}^m$ such that
\begin{equation}\label{fh}
   f(t)=-D^\top \bar h(t)
\end{equation}
is our another assumption. As we further clarify in Remark~\ref{remrem}, the proof of Theorem~\ref{moreauthm} implies that assumption (\ref{fh}) is equivalent to 
\begin{equation}\label{balance!}
   f_1(t)+...+f_n(t)=0.
\end{equation}

\noindent Introducing 
\begin{equation}\label{U}   U=\left\{x\in D\mathbb{R}^n:R^\top x=0\right\},\qquad  V=A^{-1}U^\perp,
\end{equation}
where 
\begin{equation}\label{Uperp}
U^\perp=\left\{y\in\mathbb{R}^m:\left<x,y\right>=0,\ x\in U\right\},
\end{equation}
the space $V$
will be the orthogonal complement of the space $U$ in the sense of the scalar product
\begin{equation}\label{product}
   (u,v)_A=\left<u,Av\right>.
\end{equation}
Therefore, any element $x\in\mathbb{R}^m$ can be uniquely decomposed 
as
$$
   x=P_U x+P_V x,
$$
where $P_U$ and $P_V$ are linear (orthogonal in sense of \eqref{product}) projection maps on $U$ and $V$ respectively.
Define
\begin{align}
  g(t)&=P_V DL l(t),\label{g(t)}\\
 h(t)&=P_U A^{-1}\bar h(t),\label{h}\\ 
  N_C^A(x)&=\left\{\begin{array}{ll}\left\{\xi\in\mathbb{R}^m:\left<\xi,A(c-x)\right>\le 0,\ {\rm for\ any }\ c\in C\right\},&\quad {\rm if}\ x\in C,\\
   \emptyset,& \quad {\rm if}\ x\not\in C,
\end{array}\right.\nonumber\\
  \Pi(t)&=A^{-1}C+h(t)-g(t),\label{Pi(t)}
\end{align}

and consider
the following differential inclusions 
\begin{align}
    -\dot y&\in N^A_{\Pi(t)\cap V}(y),\label{sw1}\\
 \dot z&\in \left(N^A_{\Pi(t)}(y)+\dot y\right)\cap U,\label{sw2}
\end{align}
with initial conditions
\begin{align}
y(0)&\in \Pi(0)\cap V,\label{eq:initialElastic}\\
z(0)&\in U.\label{initial}
\end{align}

\noindent The function $g(t)$ will be termed the {\it effective displacement-controlled loading.} Similarly, $h(t)$ is termed the {\it effective stress loading.} 

\vskip0.2cm

\noindent In what follows we are going to establish an equivalence between systems (\ref{eq1})-(\ref{eq5}) and (\ref{sw1})-(\ref{initial}).

\vskip0.2cm

\noindent
According to Moreau \cite[Proposition of \S6.d]{moreau}, the problem (\ref{sw2}), (\ref{initial}) admits an absolutely continuous (possibly non-unique) solution $z$ on $[0,T]$ for any absolutely continuous solution $y$ of \eqref{sw1},     
   \eqref{eq:initialElastic} defined on $[0,T]$. The analysis of the dynamics of the elastic deformation $e(t)$ therefore reduces to the analysis of the solution $y$ of the sweeping process (\ref{sw1}). In particular, stabilization of (\ref{sw1}) will imply stabilization of both  elastic deformations $e(t)=(e_1(t),...,e_m(t))^\top$ and stresses $s(t)=(s_1(t),...,s_m(t))^\top$ of springs.

\begin{theorem} \label{moreauthm} Let $D$ be the kinematic matrix of a connected network of $m$ elastoplastic springs on $n$ nodes. Let $R$ be a matrix of  incidence vectors of $q$ displacement-controlled constraints, which are independent in the sense of (\ref{A2}). 
Assume that the stress loading doesn't exceed the safe load bounds, i.e.  
\begin{equation}\label{A4}
\left(C+Ah(t)\right)\cap U^\perp\not=\emptyset,\ { for \ all\ }t\in[0,T],
\end{equation}
holds for $C$, $U$ and $h$ as defined in (\ref{formula2}), (\ref{U}), and (\ref{h}). If $(s(t),e(t),p(t),r(t))$ is a solution of the variational system (\ref{eq1})-(\ref{eq5}) on $[0,T],$ then 
\begin{equation}\label{changev}
\begin{array}{l}
  y(t)=e(t)+h(t)-g(t),\\
   z(t)=e(t)+p(t)+h(t)-g(t)
\end{array} 
\end{equation}
is a solution of the sweeping process (\ref{sw1})-(\ref{initial}) on $[0,T].$ Conversely, if $(y(t),z(t))$ is a solution of (\ref{sw1})-(\ref{initial}) then  $(e(t),p(t))$ found from (\ref{changev}) is a solution of (\ref{eq1})-(\ref{eq5})  with $s(t)=Ae(t)$ and with some suitable $r(t).$
\end{theorem}



\noindent We refer the reader to Had-Reddy \cite{han} for formulations of the safe load condition (\ref{A4}) in the context of classical (continuum) theory of plasticity.

\begin{remark}\label{A4A4}
Since $P_U g(t)=0$, condition (\ref{A4}) is equivalent to assuming 
$\Pi(t)\cap V \not=\emptyset,$ $t\in[0,T].$
\end{remark}

\begin{figure}[h]
\begin{center}
\includegraphics[scale=0.5]{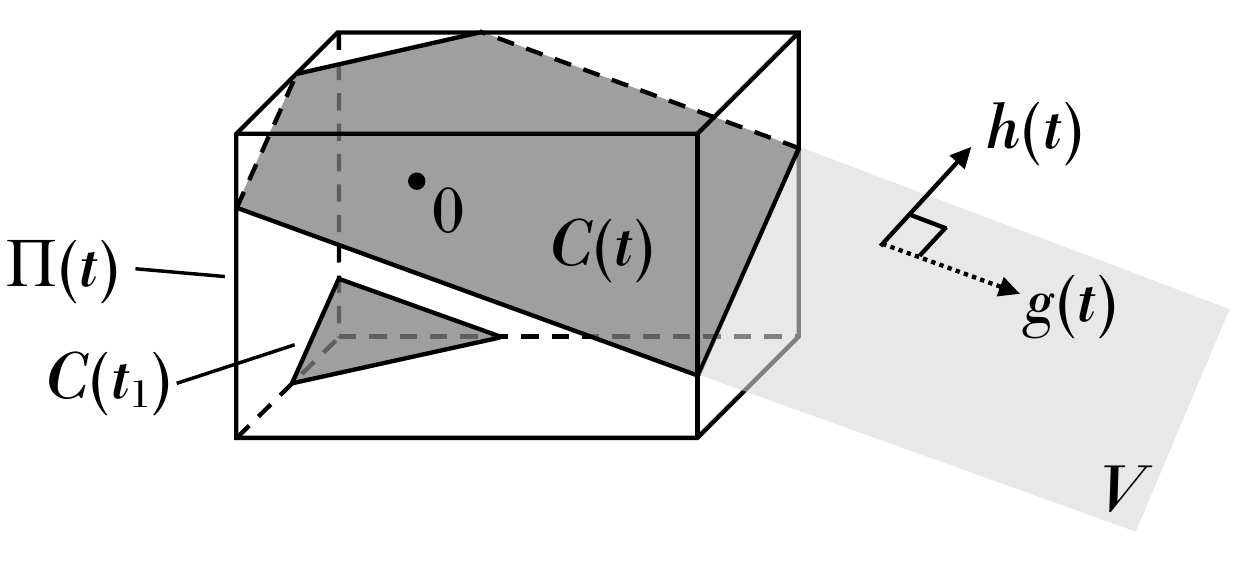} 
\end{center}
\vskip-0.3cm
\caption{Moving constraint   for different values of time. \hfill {\color{white}a} \protect\linebreak 
   }
\label{fff} 
\end{figure}

\begin{remark}
\noindent Condition (\ref{A4}) always holds when $h(t)\equiv 0$ because $0\in C$ and $0\in U^\perp$. Geometrically, condition (\ref{A4}) means that the parallelepiped $\Pi(t)$ and the hyperplane $V$ in Fig.~\ref{fff} do intersect. Mechanically, condition (\ref{A4})  accounts for the fact that the stresses of the elastoplastic springs are bounded and cannot balance arbitrary large stress loadings.
\end{remark}

\begin{remark}
\label{rem:uniqueness} The function $g(t)$ and the matrix $P_VDL$ don't depend on the choice of matrix $L$. Indeed, let $\tilde g(t)$ be the function $g(t)$ obtained by replacing  $L$ by $\tilde L$. 
Then using \eqref{U} we get
$$
0_{q\times q}=R^\top D(L-\tilde L)=R^\top (P_U D(L-\tilde L)+P_V D(L-\tilde L))
=R^\top P_V D(L-\tilde L),
$$
so $P_V D(L-\tilde L)\mathbb{R}^q\subset U$. Therefore, $P_V D(L-\tilde L)=0_{m\times m}$ and $P_V DL=P_V D\tilde L$. The conclusion about $g(t)$ follows from \eqref{g(t)}.
\end{remark}

In Appendix~\ref{appendB} we offer a diagram (fig.~\ref{fig:StructureOfSpace}) showing graphically how the spaces $U,V$ and the moving set $\Pi(t)\cap V$ are constructed.

\begin{figure}[h]\center
\includegraphics[scale=0.65]{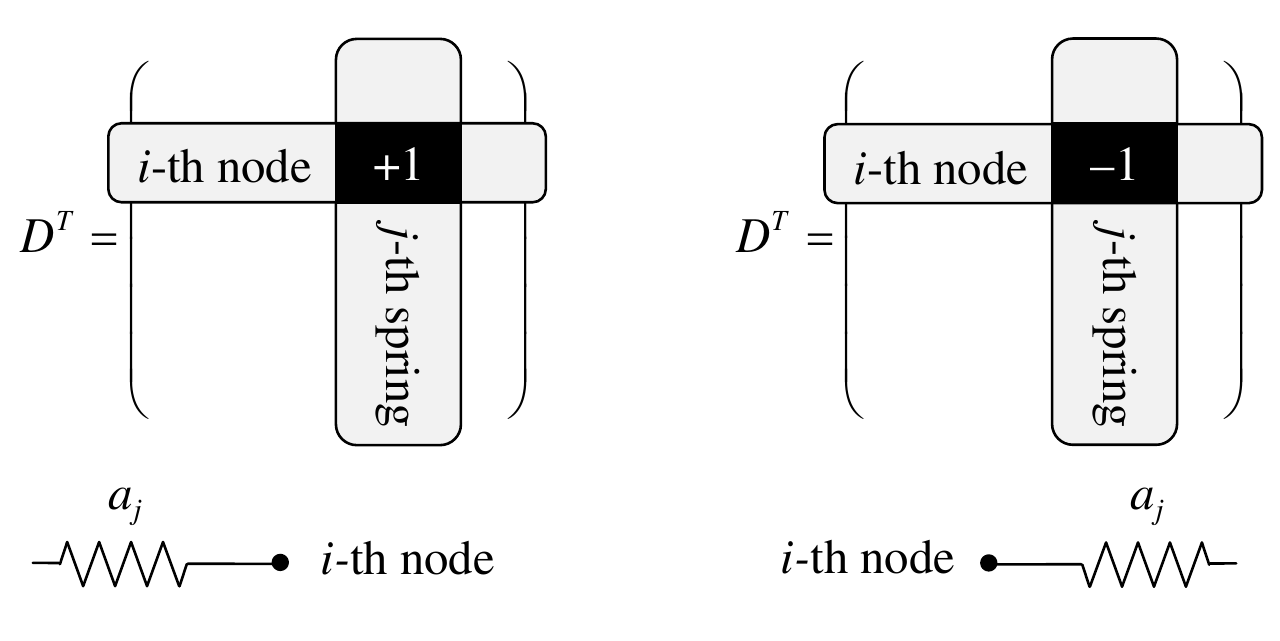}
\caption{\footnotesize The meaning of the columns and rows of matrix $D^\top.$ The cell equals $+1$, if the $i$-th node is the right endpoint for spring $j$. Conversely, the cell equals $-1$, if the $i$-th node is the left endpoint for spring $j$.} \label{fig6}
\end{figure}

\noindent {\bf Proof of Theorem~\ref{moreauthm}.}  The system of (\ref{D}) and (\ref{G}) is equivalent to 
\begin{equation}\label{rr}
   e(t)+p(t)\in U^l(t),\quad{\rm where}\quad U^l(t)=\left\{x\in D\mathbb{R}^n:R^\top x=l(t)\right\}. 
\end{equation}
Applying the both sides of (\ref{barxi}) to $l(t)$, we get $R^\top DL l(t)=l(t)$, which implies $DL l(t)\in U^l (t).$ Therefore, 
$$
  U^l(t)=U+DLl(t)
$$
and (\ref{rr}) can be rewritten as
$$
   e(t)+p(t)\in U+DL  l(t),
$$
or, equivalently,
\begin{equation}\label{ach1}
  e+p\in U+g(t).
\end{equation}
\noindent By the definition of matrix $D$, the $i$-th line has $+1$ ($-1$) at those nodes which are right (left) endpoints for the $j$-th spring, see the illustration at fig.~\ref{fig6}.
Therefore,
\begin{equation}\label{DTs}
  s_1s^1+\ldots+s_ms^m=-D^\top s.
\end{equation}
\vskip-0.4cm

\vskip0.2cm

\noindent Now we are going to demonstrate that it is enough to use the description of displacement-controlled loadings $R^k$(which is it terms of springs lengths), i.e
\begin{equation}\label{weclaim-r-R}
r_1r^1+\ldots+r_qr^q=-D^\top Rr,\qquad {\rm where}\quad r=\left(r_1,\ldots,r_q\right)^\top.
\end{equation}
Indeed, by Fig.~\ref{fig6}, the matrix  $\left(\begin{array}{c} -D\\ (r^k)^\top\end{array}\right)^\top$ is the incidence matrix of the oriented graph of springs $s_1,...,s_m$ on nodes $1,...,n$ supplemented with a virtual spring connecting the nodes $\xi_{I_k}<\xi_{J_k}$ (see section \ref{ssect:DPL} earlier).
We can now use this virtual spring in order to close the chain of springs given by the incidence vector $R^k$ and obtain a directed cycle where the direction from the node $I_k$ to the node $J_k$ disagree with the direction of the virtual spring. The incidence vector of this cycle is 
$\left(\begin{array}{c} R^k\\ -1\end{array}\right).$ According to \cite[p.~57]{bapat}, we now have 
$$\left(\begin{array}{c} -D\\ (r^k)^\top\end{array}\right)^\top\left(\begin{array}{c} R^k\\ -1\end{array}\right)=0,\qquad k\in\overline{1,q},$$
from which \eqref{weclaim-r-R} follows.

\vskip0.2cm

\noindent Therefore, taking into account (\ref{DTs}), one concludes that (\ref{eq5}) can be rewritten as 
\begin{equation}\label{Ivan11}
   -D^\top s-D^\top R r+f(t)=0,
\end{equation}
which has a solution $(s(t),r(t))$ if and only if  
$$
   s(t)+Rr(t)+\bar h(t)\in {\rm Ker}~D^\top.
$$
Keeping $s(t)$ fixed, the latter inclusion can be solved for $r(t)\in\mathbb{R}^q$ if and only if 
\begin{multline}\label{balance}
   s(t)+\bar h(t)\in {\rm Ker}~D^\top + R\,\mathbb{R}^q=\left(D\mathbb{R}^n\right)^\perp +\left\{x\in\mathbb{R}^m:R^\top x=0\right\}^\perp=\\
=\left(D\mathbb{R}^n\cap\left\{x\in\mathbb{R}^m:R^\top x=0\right\}\right)^\perp=U^\perp,
\end{multline}
see e.g. Friedberg et al \cite[Exercise~17, p.~367]{algebrabook} for the property 
${\rm Ker}~D^\top=\left(D\mathbb{R}^n\right)^\perp$.
If $s(t)$ satisfies (\ref{balance}), then by (\ref{h})
\begin{equation}\label{balance1}   
\begin{array}{r} s(t)+Ah(t)=s+A\left(P_U A^{-1}\bar h(t)+P_V A^{-1}\bar h(t)-P_V A^{-1}\bar h(t)\right)\in \\
\in s+AA^{-1}\bar h(t)+AV=s+\bar h(t)+U^\perp\in U^\perp.
\end{array}
\end{equation}
Vice versa, if $s(t)$ satisfies (\ref{balance1}) then
$$
\begin{array}{r}
  s(t)+\bar h(t)=s+A\left(P_UA^{-1}\bar h(t)+P_V A^{-1}\bar h(t)\right)=\\
=s(t)+Ah(t)+AP_VA^{-1}\bar h(t)\in U^\perp,
\end{array}
$$
which is (\ref{balance}).
By applying $A^{-1}$ to (\ref{balance1}), we get
\begin{equation}\label{ach2}
  e+h(t)\in V.
\end{equation}
Since $g(t)\in V$ and $h(t)\in U$ we can rewrite (\ref{ach1}), (\ref{ach2}) and (\ref{formula2}) as
\begin{eqnarray*}
e+p-g(t)+ h(t)&\in&  U,\\
e+h(t)-g(t)&\in&V,\\
\dot p&\in& N_C(Ae).
\end{eqnarray*}
Introducing the change of the variables  (\ref{changev})
we have $p=z-y$ and using the substitution $e=y-h(t)+g(t)$ 
\begin{equation}\label{temp}
\begin{array}{l}
\dot z-\dot y\in N^A_{A^{-1}C}(y-h(t)+g(t)),\\
 z \in  U,\\
 y\in  V.
\end{array}
\end{equation}
Let $(y,z)$ be a solution of (\ref{temp}). Since $z\in U$ we have $-\dot z\in -U=U=V^{\perp_A}=N^A_V(y)$, where $V^{\perp_A}$ is the orthogonal complement of $V$ in the sense of the scalar product $(\cdot,\cdot)_A,$ and the inclusion  (\ref{sw1}) computes as follows:
$$
\begin{array}{r}
-\dot y\in N^A_{A^{-1}C}(y-h(t)+g(t))-\dot z\in N^A_{A^{-1}C+h(t)-g(t)}(y)+N^A_V(y)=\\
=N^A_{(A^{-1}C+h(t)-g(t))\cap V}(y),
\end{array}
$$
where the last equality holds due to (\ref{A4}) (where both intersecting sets are polyhedral, we use \cite[Corollary 23.8.1]{Rockafellar}). 
The inclusion (\ref{sw2}) follows by combining  
$$
  \dot z=\dot z-\dot y+\dot y\in N^A_{A^{-1}C+h(t)-g(t)}(y)+\dot y,
$$
with the property $z(t)\in U$ observed in (\ref{temp}).  
Vice versa, if $(y,z)$ is a solution of (\ref{sw1})-(\ref{sw2}), then
\begin{eqnarray*}
&& -y\in V,\\
&& \dot z\in N^A_{A^{-1}C+h(t)-g(t)}(y)+\dot y,\\
&& \dot z\in U,
\end{eqnarray*}
which implies (\ref{temp}) when combined with (\ref{initial}).
\qed

\begin{remark}\label{remrem} Having the proof of Theorem~\ref{moreauthm} behind, we can now clarify why the $n$ equations (\ref{eq5}) of static balance in nodes is equivalent to just one equation (\ref{balance!}). Indeed, as it follows from the proof of Theorem~\ref{moreauthm}, equation (\ref{eq5}) is equivalent to (\ref{Ivan11}) which  necessarily means that
\begin{equation}\label{newequiv} f(t)\in D^\top\mathbb{R}^m.
\end{equation}
It remains to show that (\ref{newequiv}) is equivalent to (\ref{balance!}).  
Since ${\rm dim\,Ker}\,D+{\rm rank}\,D=n$, by rank-nullity theorem (see e.g.  \cite[Theorem~2.3]{algebrabook}) and ${\rm rank}\, D=n-1$ by Bapat \cite[Lemma~2.2]{bapat}  (the rank of the incidence matrix of a connected graph is one less the number of nodes), one has $\dim {\rm Ker}\,D =1.$ On the other hand,  $D(1,...,1)^\top=0$
by inspection. Therefore,  $D^\top\mathbb{R}^m=({\rm Ker} \hskip0.05cm  D)^\perp=\left\{x:(1,...,1)x=0\right\},$ i.e. 
(\ref{newequiv}) is equivalent to (\ref{balance!}).  
\end{remark}

\noindent We acknowledge that the ideas of the proof of Theorem~\ref{moreauthm} are due to Moreau \cite{moreau}, who however worked in abstract configuration spaces and didn't give details that relate the sweeping process (\ref{sw1}) to  networks of connected springs (\ref{eq1})-(\ref{eq5}).

\vskip0.2cm

\noindent Formulas (\ref{g(t)})-(\ref{Pi(t)}) establish a connection between mechanical properties of  applied loading and geometric properties of the moving constraint $\Pi(t)\cap V$. Specifically, varying the stress loading $f(t)$ moves $\Pi(t)$ in the direction perpendicular to $V$ in the sense of the scalar product (\ref{product}). In contrast, varying the displacement-controlled loading $l(t)$ moves $\Pi(t)$ in the direction parallel $V.$ We also see that the variety of possible perpendicular motions coming from $f(t)$ is limited by the dimension of the space $U,$ which will be computed in section~\ref{sec31} (Lemma~\ref{n-q-1}). The dimension of possible directions for the parallel motion in $V$ is not always $\dim V,$ but is related to the rank of matrix $\bar L$, which we compute in section \ref{sec33}, see formula (\ref{g(t)ex}). 

\subsection{Sweeping processes of particular elastoplastic systems}\label{sec31}

In this section we consider particular networks of elastoplastic springs and offer a guideline that can be used to derive the associated sweeping process (\ref{sw1}) in closed form. 
\vskip0.2cm

\noindent The following lemma will be used to compute the dimension of $U.$

\begin{lemma}\label{n-q-1} If (\ref{A2}) is satisfied, then \begin{equation}\label{dimU}
\dim U=n-q-1.
\end{equation}
\end{lemma}

\noindent {\bf Proof.} Let $E=D\mathbb{R}^n.$ Viewing $R^\top$ as a linear map from $E$ to $\mathbb{R}^q$ the rank-nullity theorem (see e.g. Friedberg et al \cite[Theorem~2.3]{algebrabook}) gives
$$
   {\rm dim\,Ker}\,R^\top+{\rm rank}\hskip0.05cm R^\top=\dim E,
$$
where ${\rm dim\,Ker}\,R^\top=\dim U$ by (\ref{U}), 
${\rm rank}\,R^\top=q$ by (\ref{A2}), and $\dim E =n-1$ by Bapat \cite[Lemma~2.2]{bapat}. \qed

\vskip0.2cm

\noindent {\bf Example.} Consider a one-dimensional network of 3 springs on 4 nodes with the kinematic matrix $D$ provided by the map
\begin{equation}\label{Dex1}
D\xi=\begin{pmatrix}
\xi_2-\xi_1\\
\xi_3-\xi_2\\
\xi_4-\xi_3
\end{pmatrix}=\begin{pmatrix}
-1 & 1 &0 &0\\
0  & -1 & 1&0\\
0&0&-1&1
\end{pmatrix}\begin{pmatrix}
\xi_1\\
\xi_2\\
\xi_3\\
\xi_4
\end{pmatrix},
\end{equation}

\noindent some $3\times 3$ diagonal matrix $A$ of Hooke's coefficients and some intervals $[c_i^-,c_i^+]$, $i\in\overline{1,3}$, of elasticity bounds. Assume that displacement-controlled loading $l(t)\in\mathbb{R}^2$ is given by the incidence vectors 
\begin{equation}\label{Rex1}
(R_1,R_2)=R=\left(\begin{array}{ll} 0 & 1\\ 1 & 1\\ 1 & 0\end{array}\right),
\end{equation}
\begin{figure}[h]\center
\vskip-0.6cm
\includegraphics{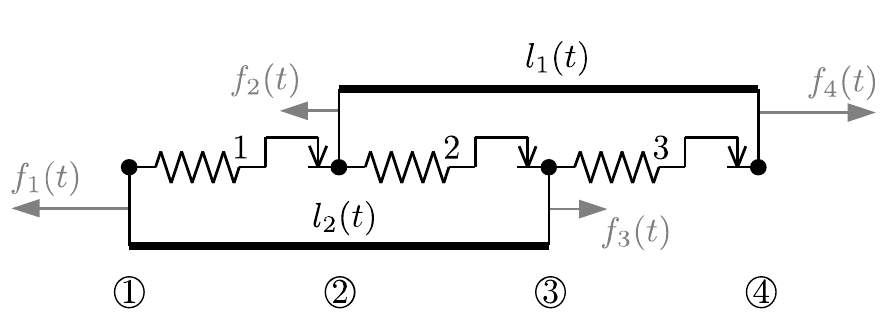}\vskip-0.0cm
\caption{\footnotesize A one-dimensional network of 3 springs on 2 nodes with 2 length locking constraints. The circled figures stays for numbers of nodes. The regular figures are the numbers of springs. The figure shows just one possible option for the directions (and magnitudes) of the forces in gray color.
} \label{ex1fig}
\end{figure}
\noindent see Fig.~\ref{ex1fig}. To examine the shapes of the associated moving set $\Pi(t)\cap V$, we find out the eligible values of the function $h(t)$. 

\vskip0.2cm

\begin{figure}[h]\center
\vskip-0.0cm
\includegraphics[scale=0.5]{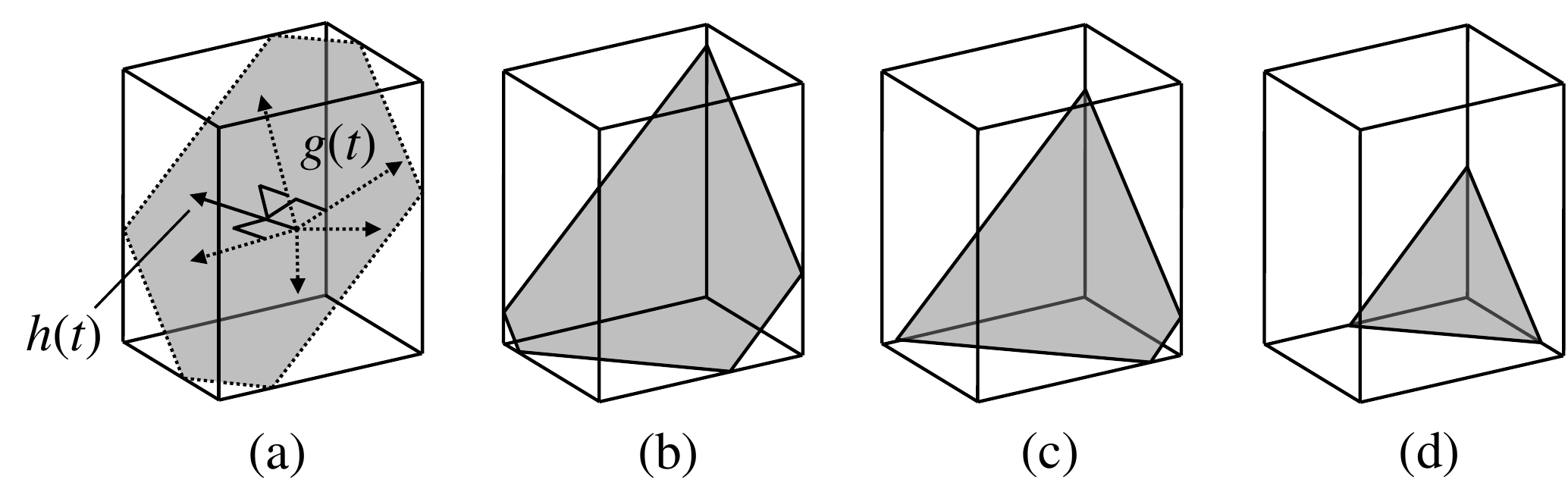}\vskip-0.0cm
\caption{\footnotesize Shapes of the moving constraint $\Pi(t)\cap V$ for the sweeping process of the network of Fig.~\ref{ex1fig} with parameters $c_1^+=-c_1^-=1, c_2^+=-c_2^-=1.3, c_3^+=-c_3^-=1.6, a_1=a_2=a_3=1$ for different values of stress loading $h(t)=(1,-1,1)^\top t$: 
\hskip5cm
a) $t=0$, b) $t=0.32,$ c) $t=0.5$, d) $t=0.8$. Figure (a) also features the possible directions of the function $g(t)$ (dotted vectors) and the possible direction of the function $h(t)$ (solid vector) that represent displacement-controlled and stress-controlled loading respectively.} \label{ex1fig_}
\end{figure}

\noindent From (\ref{h}) we conclude that eligible stress-controlled loading $f(t)$ lead to $h(t)$ given by 
\begin{equation}\label{h(t)}
  h(t)=U_{basis}H(t),
\end{equation}
where $~U_{basis}$ is the $m\times{\dim U}-$matrix of the vectors of a basis of $U$ and $H:[0,T] ~\to\mathbb{R}^{\dim U}$ is any absolutely continuous function.By \eqref{dimU} and \eqref{U}, there should exist an $n\times (n-q-1)-$matrix $M$ such that 
\begin{equation}\label{RTDM}
  R^\top DM=0 \quad {\rm and}\quad {\rm rank}(DM)=n-q-1
\end{equation}
which allows to introduce $U_{basis}$ as
\begin{equation}\label{Ubasis}
  U_{basis}=DM.
\end{equation}
Getting back to the matrices $D$ and $R$ given by (\ref{Dex1}) and (\ref{Rex1}) one has $\dim U=n-q-1=4-2-1=1$. A possible $4\times 1-$matrix that solves (\ref{RTDM}), the respective $U_{basis}$ found from (\ref{Ubasis}), and the respective function $h(t)$ given by (\ref{h(t)}) are then read as
\begin{equation}\label{UbasisEx1}
   M=\left(\begin{array}{c} 0 \\ 1 \\ 0 \\ 1\end{array}\right),\quad U_{basis}=\left(\begin{array}{c} 1\\ -1 \\ 1\end{array}\right),\quad h(t)=\left(\begin{array}{c} 1\\ -1 \\ 1\end{array}\right)H(t),
\end{equation}
where $H$ is an arbitrary absolutely continuous function from $[0,T]$ to $\mathbb{R}.$ Fig.~\ref{ex1fig_} illustrates the shapes of $\Pi(t)\cap V$ for different constant values of $H(t),$ where according to (\ref{U}) we considered
\begin{equation}\label{Vort}
  V={\rm Ker}\, \left(U_{basis}^\top A\right)=\left(\begin{array}{c}a_1\\ -a_2\\ a_3\end{array}\right)^\perp.
\end{equation}

\subsection{Bounds on the stress loading to satisfy the safe load condition}\label{sec32}
In this section we are dealing with a general elastoplastic system again.
To verify condition (\ref{A2}) it is sufficient to check that displacement-controlled lengths $l_i(t)$ can be varied independently one from another. Computational algorithms to verify safe load condition (\ref{A4}) for particular systems is a standard topic of computational geometry, see e.g. Bremner et al \cite{bremner}. In this section we derive analytic conditions which allow to spot classes of elastoplastic systems for which the safe load condition holds.

\begin{proposition}\label{cor1} In order for the safe load condition (\ref{A4}) of Theorem~\ref{thm1} to hold for some $t\ge 0,$ it is sufficient to assume that 
\begin{equation}\label{A4_}
   -Ah(t)\in C.
\end{equation}
\end{proposition}

\noindent {\bf Proof.} 
\noindent In order to show that (\ref{A4_}) implies (\ref{A4}), it  is sufficient to observe that $0\in U^\perp$ and that (\ref{A4_}) yields $0\in C+Ah(t)$.\qed

\vskip0.2cm

\begin{definition}\label{def2}
We will say that a {\it spring $i$ is blocked by displacement-controlled loadings}, if the family of displacement-controlled loadings $\{l_j\}_{j=1}^q$ contains a chain that connects one end of spring $i$ with its other end. 
\end{definition}

\begin{lemma}\label{lemu0} Assume that in a given elastoplastic system the number $q$ of displacement-controlled loadings is 2 less the number of nodes. If none of the springs of the elastoplastic system (\ref{eq1})-(\ref{eq5}) is blocked by displacement-controlled loadings, then 
$$
   x_i\not= 0 \mbox{ for any }i\in\overline{1,m},\ x\in U\backslash\{0\}.
$$
\end{lemma}

\noindent {\bf Proof.} Recall, that $I_k$ and $J_k$ are the left and right endpoints respectively of the displacement-controlled constraint $l_k(t).$ Consider the matrix $D_1$ obtained from matrix $D$ by combining the column $I_k$ and the column $J_k$ as follows: 1) add the values of column $J_k$ to the respective values of column $I_k$, 2) delete the column $J_k.$ Then,
$$
   \left\{D_1\xi:\xi\in\mathbb{R}^{n-1}\right\}=\left\{D\xi:(R^k)^\top D\xi=0,\ \xi\in\mathbb{R}^n\right\}.
$$
Moreover, each row of $D$ has exactly one element $1$ and $-1$ and none of the springs are blocked therefore at least one of each pair of summands at step 1) is zero. Matrix $D_1$ is the kinematic matrix for a new elastoplastic system that is obtained from elastoplastic system (\ref{eq1})-(\ref{eq5}) by merging the nodes $I_k$ and $J_k$ together and, thus, by reducing the number of nodes by 1. Accordingly, the new elastoplastic system features only $q-1$ displacement-controlled loadings and the indices $\{I_i,J_i\}_{i=1}^{q-1}$ are now from $\overline{1,n-1}.$ 

\vskip0.2cm

\noindent Repeating this process through all the incidence vectors $\{R^k\}_{k=1}^q,$ where $q=n-2$ by Lemma \ref{n-q-1}, we obtain 
$$
U=\left\{D\xi:R^\top D\xi=0,\ \xi\in\mathbb{R}^n\right\}=\left\{\bar D\xi:\xi\in\mathbb{R}^{n-q}\right\}=\left\{\bar D\xi:\xi\in\mathbb{R}^2\right\},
$$
where $\bar D$ is the kinematic matrix of the reduced elastoplastic system that is obtained from the original one by merging node $I_k$ with node $J_k$ trough $k\in\overline{1,q}.$ 

\vskip0.2cm

\noindent Since the reduced elastoplastic system has only two nodes ($q=n-2$), all the displacement-controlled constraints of the original system split into at most two connected components, which shrink into these two nodes under the proposed  reduction process. If spring $i$ is not blocked by displacement-controlled loadings, then the endpoints of spring $i$ belong to different connected components introduced. Therefore, the endpoints of spring $i$ are two different nodes of the reduced elastoplastic system, which implies
\[
   u_i\not=0,\  i\in\overline{1,m},\quad \mbox{for any } u=\bar D\xi\ \text{such that } \xi\in\mathbb{R}^2, \xi_1\neq\xi_2.
\]
 If $\xi_1=\xi_2$ then $u=\bar D \xi=0$. The proof of the lemma is complete.\qed

\begin{proposition}\label{cor2} Assume that the conditions of Lemma \ref{lemu0} hold. Let 
$\bar u$ be an arbitrary nonzero fixed vector of $U$ {\rm(}$\dim U=1$ by Lemma~\ref{n-q-1}{\rm)} and consider
$$
  \bar c^+=\left(\begin{array}{c}
c_1^{\sign(\bar u_1)}\\
\vdots\\
c_m^{\sign(\bar u_m)}
\end{array}\right),\qquad \bar c^-=\left(\begin{array}{c}
c_1^{-\sign(\bar u_1)}\\
\vdots\\
c_m^{-\sign(\bar u_m)}
\end{array}\right),
$$
where $c_i^{-1}$ denotes $c_i^-$ and $c_i^{+1}$ denotes $c_i^+$.
Then, for each fixed $t\ge 0,$ the safe load condition (\ref{A4}) holds if and only if 
\begin{equation}\label{A4__}
   \left<\bar u,\bar c^+ +Ah(t)\right>\cdot \left<\bar u,\bar c^-+Ah(t)\right>\le 0.
\end{equation}
\end{proposition}

\vskip0.2cm

\noindent {\bf Proof.} We first show that (\ref{A4__}) implies (\ref{A4}). Assume that (\ref{A4}) doesn't hold for some $t\in[0,T].$ Therefore, by convexity of $C$, either $\left<\bar u,x+Ah(t)\right>>0$ for all $x\in C$ or $\left<\bar u,x+Ah(t)\right><0$ for all $x\in C$. In either case we conclude $\left<\bar u,\bar c^++Ah(t)\right>\cdot \left<\bar u,\bar c^-+Ah(t)\right>>0$ because $\bar c^+,\bar c^-\in C,$ which contradicts (\ref{A4__}).

\vskip0.2cm

\noindent Let us now show that (\ref{A4}) implies (\ref{A4__}). Indeed, since
$$
   \bar u_i \bar c_i^-\le \bar u_i c_i^j\le \bar u_i \bar c_i^+,\quad \mbox{for any }i\in\overline{1,m},\ j\in\{-1,+1\},
$$
we have
$$
   \left<\bar u,\bar c^-+Ah(t)\right>\le \left<\bar u,x+Ah(t)\right>\le \left<\bar u,\bar c^++Ah(t)\right>,\quad \mbox{for any }x\in C.
$$
The latter inequality takes the required form (\ref{A4__}) when one plugs $x$ satisfying 
$\left<\bar u,x+Ah(t)\right>=0$, which exists because of (\ref{A4}).\qed

\begin{remark}\label{remark4} Considering the left-hand-side of (\ref{A4__}) as a polynomial $P\left(\left<\bar u,Ah(t)\right>\right)$ in $\left<\bar u,Ah(t)\right>$, we see that the branches of the polynomial are pointing upwards. Therefore, condition (\ref{A4__}) is the requirement for $\left<\bar u,Ah(t)\right>$ to stay strictly between the roots of the polynomial. The roots of $P\left(\left<\bar u,Ah(t)\right>\right)$ are given by $\left<\bar u,Ah(t)\right>=-\left<\bar u,\bar c^-\right>>0$ and $\left<\bar u,Ah(t)\right>=-\left<\bar u,\bar c^+\right><0$. Therefore, (\ref{A4__}) is equivalent to
$$
  -\left<\bar u,\bar c^+\right>\le\left<\bar u,Ah(t)\right>\le-\left<\bar u,\bar  c^-\right>,
$$
which highlights that (\ref{A4__}) is a restriction on the magnitude of $\left<\bar u,Ah(t)\right>.$
\end{remark}

\noindent Proposition~\ref{cor2} can be e.g. applied to the one-dimensional network of Fig.~\ref{ex1fig}, where $\dim U=1$ as we noticed earlier.

\vskip0.2cm

\noindent {\bf Example (continued).} For the elastoplastic system of Fig.~\ref{ex1fig} one can consider $\bar u=(1,-1,1)^\top$ and using 
(\ref{UbasisEx1}) obtain
\begin{equation}\label{ex1form}
  \bar c^+=\left(\begin{array}{c}
c_1^+\\ c_2^- \\ c_3^+\end{array}\right),\quad 
\bar c^-=\left(\begin{array}{c}
c_1^-\\ c_2^+ \\ c_3^-\end{array}\right),\quad 
\left<\bar u,Ah(t)\right>=(a_1+a_2+a_3)H(t).
\end{equation}
Based on Remark~\ref{remark4} the necessary and sufficient condition for safe load condition (\ref{A4__}) to hold is then
\begin{equation}\label{ex1_safeload}
    -c_1^++c_2^--c_3^+\le (a_1+a_2+a_3)H(t) \le  -c_1^-+c_2^+-c_3^-.
\end{equation}

\subsection{Condition on the displacement-controlled loading to eliminate constant solutions} \label{sec33}

\noindent Next proposition gives conditions to ensure that any point $x$ which belongs to the moving set $\Pi(t)\cap V$ of sweeping process (\ref{sw1}) at some initial time $t=t_1$ will lie outside $\Pi(t)\cap V$ at time $t=t_2.$ These conditions will, therefore, rule out the existence of constant solutions.

\begin{proposition}\label{prop3} Assume that conditions of Theorem~\ref{moreauthm} hold. If 
\begin{equation}\label{formulaprop3}
\|A^{-1}c^--A^{-1}c^+\|_A<\left\|g(t_1)-g(t_2)\right\|_A,
\end{equation}
for some $0\le t_1< t_2$, 
where 
$$
   \|x\|_A=\sqrt{\left<x,Ax\right>},\quad c^-=(c_1^-,...,c_m^-)^\top,\quad c^+=(c_1^+,...,c_m^+)^\top,
$$
then sweeping process (\ref{sw1}) doesn't have any solutions that are constant on $[t_1,t_2].$
\end{proposition}

\noindent {\bf Proof.} The claim follows by showing that 

$$
\left(\Pi(t_1)\cap V\right)\cap \left(\Pi(t_2)\cap V\right)=\emptyset
$$

Since $h(t)\in U$ we have 
$$
 \Pi(t)\cup V=\left(A^{-1}C+h(t)-g(t)\right)\cap V\subset P_V A^{-1}C-g(t),\quad t\in[t_1,t_2],
$$
 and it is sufficient to prove that the sets
$$
P_V A^{-1}C-g(t_1)\ \
{\rm and}\ \ 
P_V A^{-1}C-g(t_2)\ \ \mbox{don't intersect}.
$$
The latter will hold, if the diameter of the set $P_V A^{-1}C$ is smaller than the distance between $g(t_1)$ and $g(t_2),$ which fact will now be established.

\vskip0.2cm

\noindent Since $P_V$ is the orthogonal projection in the sense of the scalar product $(x,y)_A=\left<x,Ay\right>,$ we have (see e.g. Conway \cite[Theorem 2.7~b)]{conway})
$$
    \|P_V x\|_A\le \|x\|_A,\quad x\in\mathbb{R}^m.
$$
Therefore, for any $c_1,c_2\in C$, 
\begin{eqnarray*}
 && \left\| P_V \left(A^{-1}c_1-A^{-1}c_2\right)\right\|_A\le
\left\|A^{-1}c_1-A^{-1}c_2\right\|_A\le \|A^{-1}c^--A^{-1}c^+\|_A<\\
&& < \left\|g(t_1)-g(t_2)\right\|_A.
\end{eqnarray*}
The proof of the proposition is complete.

\vskip0.2cm

\begin{remark}\label{remarkprop3} Note, the left-hand-side in the squared inequality \eqref{formulaprop3} from the statement of Proposition~\ref{prop3} can be computed as \newline
$\|A^{-1}c^--A^{-1}c^+\|_A^2=
 \left<c^--c^+,A^{-1}(c^--c^+)\right>.$
\end{remark}

\noindent In what follows we show which kind of computations is required to verify the condition of Proposition~\ref{prop3} in practice.

\vskip0.2cm

\noindent {\bf Example (continued).}  
\noindent Given the elastoplastic system of Fig.~\ref{ex1fig}, our goal is to compute the effective displacement-controlled loading $g(t)$ of \eqref{g(t)}.  

\vskip0.2cm

\noindent For the term $P_V DL$ of \eqref{g(t)} observe, that there exists 
a $\dim V\times q$-matrix
 $\bar L$ such  that 
\begin{equation}
V_{basis}\bar L=P_V DL,
\label{eq:barL0}
\end{equation}
where $V_{basis}$ is the $m\times\dim V-$matrix of the vectors of a basis of $V$. The $i$-th column of matrix $\bar L$ is the vector of the coordinates of the respective vector $P_V DL^i \in V$ in the basis $V_{basis},$ where  $L^i$ stays for the $i$-th column on matrix $L.$ Formula  (\ref{g(t)}) can therefore be rewritten as
\begin{equation}\label{g(t)ex}
   g(t)=V_{basis}\bar L l(t).
\end{equation}
Computing the effective displacement-controlled loading $g(t)$ has hereby  been turned into computing $V_{basis}$ and $\bar L$.

\vskip0.2cm

\noindent  By (\ref{dimU}),
\begin{equation}\label{dimV}
    \dim V=m-n+q+1,
\end{equation}
and according to (\ref{Vort}), $V_{basis}$ is an arbitrary matrix of $\dim V$ linearly independent columns that solves
\begin{equation}\label{Vbasis}
  (U_{basis})^\top AV_{basis}=0.
\end{equation}
For the particular matrices (\ref{Dex1}), using the earlier computed $U_{basis}$, see (\ref{UbasisEx1}), one gets $\dim V=2$, $(U_{basis})^\top=(1\ -1\ 1)$, and a possible solution to (\ref{Vbasis}) is 
$$
   V_{basis}=\left(\begin{array}{cc} 1/a_1 & 0 \\ 1/a_2 & 1/a_2 \\ 0 & 1/a_3\end{array}\right).
$$

\noindent To find $\bar L$, we observe that by (\ref{barxi}), for any $\xi\in\mathbb{R}^q$, we have 
$$
P_V DL\xi=DL\xi-P_U DL\xi \in D\mathbb{R}^n,
$$
as $P_U DL\xi\in D\mathbb{R}^n$ by definition of $U.$ Combining this relation with (\ref{barxi}) and (\ref{eq:barL0}) one gets the following equations for $\bar L$:
\begin{align}
    R^\top V_{basis}\bar L&=I_{q\times q},\label{barL1}\\
    V_{basis}\bar L\,\mathbb{R}^q&\subset D\mathbb{R}^n,\label{barL2}
\end{align}
from which $\bar L$ can be found. In Appendix~\ref{appendB} we offer a diagram (fig.~\ref{fig:StructureOfSpace}) showing the construction of $V_{basis}\bar L$ graphically.
\vskip0.2cm

For specific matrix $D$ given by (\ref{Dex1}), one has 
$D\mathbb{R}^4=\mathbb{R}^3$ (i.e. there is no geometric constraint coming from the graph of springs in this case) and so \eqref{barL2} holds for any matrix $\bar L.$ The matrix $\bar L$ is therefore a 
$2\times 2-$matrix that solves (\ref{barL1}), which has a unique solution
$$
  \bar L=\left(\begin{array}{cc} 1/a_2 & 1/a_2+1/a_3\\
1/a_1+1/a_2 & 1/a_2\end{array}\right)^{-1}.
$$
Formula (\ref{g(t)ex}), in particular, implies that, for the network of springs of Fig.~\ref{ex1fig} (where $\dim V=q=2$), the displacement-controlled constraints are capable to execute any desired motion of $C(t)$ in $V.$ 

\vskip0.2cm

\noindent Applying Proposition~\ref{prop3} and Remark~\ref{remarkprop3}, we obtain the following condition for non-existence of constant solutions. The elastoplastic system of Fig.~\ref{ex1fig} doesn't have constant solutions on $[0,T]$, if there exist $t_1,t_2\in[0,T]$ such that
\begin{equation}\label{ex1_nonconstant}
\sum_{i=1}^4 \dfrac{1}{a_i}(c_i^+-c_i^-)^2<  \left\|V_{basis}\bar L\left(l(t_1)-l(t_2)\right)\right\|^2_A.
\end{equation}

\subsection{A polyhedral description of moving sets for elastoplastic systems and reduction to subspace $V$}\label{polyhedronsection}

To give a deeper look into the possible dynamics of sweeping processes of elastoplastic systems we now rewrite moving set $\Pi(t)\cap V$ and process \eqref{sw2} in a slightly different form which is more suitable for further analysis. From 
$$
   A^{-1}C=\left\{x\in\mathbb{R}^m:\dfrac{1}{a_i}c_i^-\le x_i \le \dfrac{1}{a_i}c_i^+\right\}
$$  
we have
$$
A^{-1}C+h(t)-g(t)=\left\{x\in\mathbb{R}^m:c_i^-+a_ih_i(t)\le \left<e_i,Ax+Ag(t)\right>\le c_i^++a_ih_i(t)\right\},
$$
where $e_i\in\mathbb{R}^m$ is the vector with 1 in the $i$-th component and zeros elsewhere. Since $g(t)\in V,$ one has
$$
   \left<e_i,A x+Ag(t)\right>=\left<P_U e_i+P_V e_i,Ax+Ag(t)\right>=\left<P_Ve_i,Ax+Ag(t)\right>,\quad x\in V,
$$
and we conclude 
\begin{equation}\label{222}
   \Pi(t)\cap V=\bigcap_{i=1}^{m}\left\{x\in V:c_i^-\hskip-0.05cm+\hskip-0.05cma_ih_i(t)\hskip-0.05cm\le\hskip-0.05cm\left<n_i,Ax\hskip-0.05cm+\hskip-0.05cmAg(t)\right>\hskip-0.05cm\le\hskip-0.05cm c_i^+\hskip-0.05cm+\hskip-0.05cma_ih_i(t)\right\},
\end{equation}
where $n_i=P_V e_i$, i.e. $n_i$ are the columns of the projection matrix $P_V$.

Moreover, since $y(t)\in V$ for all $t$ and $y'(t)\in V$ for a.a. $t$ we can restrict the normal cone from \eqref{sw2} to the normal cone defined within the subspace $V$, which also appears to be the intersection of the original normal cone with $V$:
\[
N^V_{\mathcal{C}}(x):=N^A_\mathcal{C}(x)\cap V=\left\{\begin{array}{ll}\left\{\xi\in V:\left<\xi,A(c-x)\right>\le 0,\ {\rm for\ any }\ c\in \mathcal{C}\right\},&\text{ if } x\in \mathcal{C},\\
   \emptyset,&\text{ if } x\not\in \mathcal{C}.
\end{array}\right.
\]
Therefore we can restrict sweeping process \eqref{sw2} to one completely defined within $V$:
\begin{equation}
-\dot y\in N^V_{\Pi(t)\cap V}(y).
\label{eq:sw3}
\end{equation}
In the following chapters we are going to analyze dynamics of the process \eqref{eq:sw3} with the moving set in form \eqref{222}.
\section{Convergence to a periodic attractor}\label{sec4}

\subsection{Convergence in the case of a moving constraint given by an intersection of translationally moving convex sets}\label{sec41}

 \noindent In this section we establish convergence properties of 
a general sweeping process
\begin{equation}\label{swgen}
   -\dot x\in N_{C(t)}^0(x),\qquad y\in E,
\end{equation}
where $E$ is a $d$-dimensional linear vector space, $C(t)\subset E$ is convex closed set for any $t$, and 
\begin{equation}\label{normal_cone_standard}    N_C^0(x)=\left\{\begin{array}{ll}\left\{\xi\in E:\left(\xi,c-x\right)_0\le 0,\ {\rm for\ any }\ c\in C\right\},&\quad {\rm if}\ x\in C,\\
   \emptyset,& \quad {\rm if}\ x\not\in C,
\end{array}\right.
\end{equation}
where $(\cdot,\cdot)_0$ is some inner product in $E.$ 
These convergence properties are then refined in section \ref{sec1a} in the context of the particular sweeping process (\ref{sw1}).  

\vskip0.2cm

\noindent A set-valued function $t\mapsto C(t)$ is called globally Lipschitz continuous, if
\begin{equation}\label{Lip}
  d_H(C(t_1),C(t_2))\le L_C|t_1-t_2|,\ \mbox{for all}\ t_1,t_2\in\mathbb{R}, \mbox{ and for some }L_C>0,
\end{equation}
where  $d_H(C_1,C_2)$ is the Hausdorff distance  between two closed sets $C_1,C_2\in E$ defined as
\begin{equation}\label{dH}
\begin{array}{rcl}
   d_H(C_1,C_2)&=&\max\left\{  
\sup\limits_{x\in C_2} {\rm dist}(x,C_1),\sup\limits_{x\in C_1} {\rm dist}(x,C_2)
\right\}
\end{array}
\end{equation}
with ${\rm dist}(x,C)=\inf\left\{
|x-c|:c\in C
\right\}.$

\vskip0.2cm

\noindent Recall, if $C(t)$ is a globally Lipschitz continuous function with nonempty closed convex values from $E$, then the solution $x(t)$ of sweeping process (\ref{swgen}) with any initial condition $x(t_0)=x_0$ is uniquely defined on $[t_0,\infty)$ in the sense that $x(t)$ is a Lipschitz continuous function that verifies (\ref{swgen}) for a.a. $t\in[t_0,\infty)$ (see e.g. Kunze and Monteiro Marques \cite{kunze}).

\vskip0.2cm

\noindent Let us use $t\mapsto X(t,x_0)$ to denote the solution of sweeping process (\ref{swgen}) that takes the value $x_0$ at time 0. In what follows, we consider the set of $T$-periodic solutions of (\ref{swgen}) 
\begin{align}
   X&=\{X(\cdot,x_0):x_0=X(T,x_0)\}\subset C(\mathbb{R},E)\\
   X(t)&=\left\{X(t,x_0):x_0=X(T,x_0)\right\}\label{X(t)}\subset E
\end{align}
and prove that, for $T$-periodic moving constraint $C(t)$, the set $X(t)$ attracts all the solutions of (\ref{swgen}). Note that the condition $X(0,x_0)=X(T,x_0)$ implies $X(0,x_0)=X(jT,x_0),$ $j\in\mathbb{N}$, when $t\mapsto C(t)$ is $T$-periodic.

\begin{definition} A set-valued function $t\mapsto Y(t)$ is a global attractor of sweeping process (\ref{swgen}), if ${\rm dist}(x(t),Y(t))\to 0$ as $t\to\infty$ for any solution $x$ of sweeping process (\ref{swgen}).
\end{definition}

\noindent Finally, we  denote by ${\rm ri}(C)$ the relative interior of a convex set $C \subset E$, see Rockafellar \cite[\S6]{Rockafellar}.

\begin{theorem}\label{thm1} Let $t\mapsto C(t)$ be a Lipschitz continuous uniformly bounded $T$-periodic set-valued function with nonempty closed convex  values from $E$. Let $t\mapsto X(t)$ be the set of $T$-periodic solutions of sweeping process (\ref{swgen}) as defined in (\ref{X(t)}).  Then, $X\subset C([0,T],E)$ is closed and convex.  If, in addition,  $C(t)$ is an intersection of closed convex  sets $C_i$ (some of them, say, first $p$ sets, may be polyhedral) that undergo just translational motions 
\begin{equation}\label{eq:finiteIntersectionC}
  C(t)=\bigcap\limits_{i=1}^k(C_i+c_i(t)),
\end{equation}
where  $c_i(t)$ are single-valued $T$-periodic Lipschitz functions such that

\begin{equation}\label{nonempty}
   \bigcap\limits_{i=1}^p(C_i+c_i(t))\cap
   \bigcap\limits_{i=p+1}^k({\rm ri}(C_i)+c_i(t))\not = \emptyset,   
   \qquad t\in[0,T],
\end{equation}
then 
\begin{equation}\label{equivV}
\dot x(t)=\dot y(t),\quad for \ any \ x,y\in X \ and\ { for\ a.a.\ }t\in[0,T],
\end{equation}
and
$X(t)$ is a global attractor of (\ref{swgen}).
\end{theorem}

\noindent The theorem, in particular, implies that $X(t)$ cannot contain non-constant solutions, if it contains at least one constant solution.

\vskip0.2cm

\noindent The proof of theorem~\ref{thm1} is split into 3 lemmas. Lemma~\ref{lemconvex} establishes the convexity of $X$ (closedness of  $X(t)$ follows from the continuous dependence of solutions of (\ref{swgen}) on the initial condition, see \cite[Corollary~1]{kunze}). Lemma~\ref{th:constDistance} proves the statement (\ref{equivV}). Finally, the global attractivity of $X(t)$ is given by Theorem~\ref{lemattract} which is an extension of a result from Krejci \cite{Krejci1996} for convex sets 
(\ref{eq:finiteIntersectionC}).
\vskip0.2cm

\noindent In what follows, $\|\cdot\|_0$ is the norm induced by the scalar product in $E$, i.e.
\begin{equation}\label{inducednorm}
\|x\|_0=\sqrt{(x,x)_0}.
\end{equation}

\begin{lemma}\label{lemconvex} Let $t\mapsto C(t)$ be a Lipschitz continuous set-valued function with nonempty closed convex values from $E$. Then, both $X(t)\subset E$ and $X\subset C(\mathbb{R},E)$ are convex. In addition,  for any $x,y\in X\subset C(\mathbb{R},E)$
\begin{equation}
\|x(t)-y(t)\|_0 \text{ is constant in } t.
\label{eq:constantDistance}
\end{equation}
\end{lemma}

\noindent {\bf Proof.} Let $x,y\in X$. Due to monotonicity of $N_{C(t)}^0(x)$ in $x$ the distance 
$t\mapsto \|x(t)-y(t)\|_0$ 
cannot increase (see e.g. \cite[Corollary~1]{kunze}). Notice, that $t\mapsto \|x(t)-y(t)\|_0$ cannot decrease, otherwise it cannot be periodic, so \eqref{eq:constantDistance} follows.

\vskip0.2cm

\noindent For any $\theta\in(0,1)$ the initial condition $\theta x(0)+(1-\theta)y(0)$ belongs  $C(0)$ by convexity of $C$. Let $x_\theta$ be the corresponding solution. Since $t\mapsto \|x(t)-x_\theta(t)\|_0$ and $t\mapsto \|x_\theta(t)-y(t)\|_0$ are  also non-increasing, then
\begin{multline*}
\|x(t)-x_\theta(t)\|_0+\|y(t)-x_\theta(t)\|_0\leqslant \|x(0)-x_\theta(0)\|_0+\|y(0)-x_\theta(0)\|_0=\\
 =\|x(0)-y(0)\|_0=\|x(t)-y(t)\|_0.
\end{multline*}
On the other hand,  the triangle inequality yields 
$$
\|x(t)-x_\theta(t)\|_0+\|y(t)-x_\theta(t)\|_0\geqslant
 \|x(t)-y(t)\|_0
$$
and we have
\begin{equation}
   \|x(t)-x_\theta(t)\|_0+\|x_\theta(t)-y(t)\|_0=\|x(t)-y(t)\|_0 \equiv{\rm const}.
\label{eq:xThetaConst}
\end{equation}
Because none of the terms  $\|x(t)-x_\theta(t)\|_0,\|x_\theta(t)-y(t)\|_0$ can increase, both of them remain constant and positive(due to the choice of $x_\theta(0)$). Moreover, by strict convexity of the inner product space (see Narici-Beckenstein \cite[Th 16.1.4 d)]{Narici}) there is $\alpha(t)>0$ such that 
\[
x_\theta(t)-y(t)=\alpha(t)(x(t)-x_\theta(t)).
\]
We solve for $x_\theta$:
\begin{equation}
x_\theta(t)=\frac{\alpha(t)}{1+\alpha(t)}x(t)+\frac{1}{1+\alpha(t)}y(t).
\label{eq:xThetaAlpha}
\end{equation}
and substitute it to the second difference in \eqref{eq:xThetaConst}: 
\begin{gather*}
\|x(t)-x_\theta(t)\|_0=\frac{1}{1+\alpha(t)}\|x(t)-y(t)\|_0,
\end{gather*}
Both distances $\|x(t)-x_\theta(t)\|_0$ and $\|x(t)-y(t)\|_0$ are constant, hence $\alpha(t)$ is constant as well, which means that $\alpha(t)=\alpha(0)$ and due to the choice of $x_\theta(0)$ expression \eqref{eq:xThetaAlpha} becomes
$$x_\theta(t)=\theta x(t)+(1-\theta)y(t).$$
This formula, in particular, implies that $x_\theta$ is $T$-periodic. The proof of convexity of $X$ is complete.

\begin{lemma}\label{th:constDistance}
Let $t\mapsto C(t)$ be a set-valued function of the form  (\ref{eq:finiteIntersectionC})-(\ref{nonempty}) with convex closed $C_i$ and Lipschitz-continuous single valued $c_i(t)$.
Let $x$ and $y$ be two solutions of sweeping process \eqref{swgen} defined on $[0,T]$ such that \eqref{eq:constantDistance} holds for them. Then for almost all $t\in[0,T]$
\begin{equation}
\dot x(t)=\dot y(t).
\end{equation}

\end{lemma}
\noindent {\bf Proof.} 
\noindent The properties 
(\ref{eq:finiteIntersectionC})-(\ref{nonempty})
imply  
(see \cite[Corollary 23.8.1]{Rockafellar})
 that
\begin{equation}
N_{C(t)}^0(x)=\sum_{i=1}^k N_{C_i+c_i(t)}^0(x),\quad \text{for all }x\in C(t)\ {\rm and\ for\ a.a.\ }t\in[0,T].
\label{eq:additivityOfCones}
\end{equation}
Let $t\in(0,T)$ be such that $\dot x(t),\dot y(t),\dot c_i(t)$, $i\in\overline{1,k}$, exist and \eqref{eq:additivityOfCones} holds.
Property \eqref{eq:additivityOfCones} allows to spot
$\dot x^t_i,$ $\dot y^t_i$, $i\in\overline{1,k}$, such that
\begin{eqnarray*}
\dot x(t)=\sum\limits_{i=1}^k \dot x^t_i,&&\ 
-\dot x^t_i \in N_{C_i+c_i(t)}^0(x(t)),\quad i\in\overline{1,k},\\
\dot y(t)=\sum\limits_{i=1}^k \dot y^t_i,&&\ 
-\dot y^t_i \in N_{C_i+c_i(t)}^0(y(t)),\quad i\in\overline{1,k}.
\end{eqnarray*}
\noindent To show that $\|\dot x(t)-\dot y(t)\|_0=0$, consider
\begin{eqnarray}
\|\dot x(t)-\dot y(t)\|^2_0&=&\left(\dot x(t)-\dot y(t),\dot x(t)-\dot y(t)\right)_0=\nonumber\\
&=&\sum_{i=1}^k \left( \dot x_i^t,\dot x(t)-\dot c_i(t)\right)_0
+\sum_{i=1}^k\left( \dot y_i^t,\dot y(t)-\dot c_i(t)\right)_0-\label{vanish1}\\
&&\quad-\sum_{i=1}^k\left( \dot x_i^t,\dot y(t)-\dot c_i(t)\right)_0-\sum_{i=1}^k\left( \dot y_i^t,\dot x(t)-\dot c_i(t)\right)_0
.\label{vanish2}
\end{eqnarray}
For the value of $t\in(0,T)$ as fixed above, we now prove that each of sums in (\ref{vanish1})-(\ref{vanish2}) vanish.  

\vskip0.2cm

\noindent {\bf Step 1.} {\it Vanishing sums in (\ref{vanish1}).} Fix $i\in\overline{1,k}.$ By the definition of normal cone, 
\begin{equation}\label{dnc}
\left(\dot x_i^t,z+c_i(t)-x(t)\right)_0\ge 0\ {\rm and} 
\ \left(\dot y_i^t,z+c_i(t)-y(t)\right)_0\ge 0\ 
\text{for all}\  z\in C_i.
\end{equation} Considering $z=x(t+h)-c_i(t+h)\in C_i$, we observe that the function
$$  
  f(h)=\left(\dot x_i^t, x(t+h)-x(t)-(c_i(t+h)-c_i(t))\right)_0
$$
is non-negative in a neighborhood of zero. Since $f(0)=0$, we conclude that $0=f'(0)=\left(\dot x_i^t,\dot x(t) -\dot c_i(t)\right)_0.$ The relation 
$\left(\dot y_i^t,\dot y(t) -\dot c_i(t)\right)_0=0$ can be proved by analogy using the second inequality of (\ref{dnc}). 
\vskip0.2cm

\noindent {\bf Step 2.} {\it Vanishing sums in (\ref{vanish2}).} We claim that 
\begin{equation}\label{weclaim}
\sum_{i=1}^k\left(\dot x_i^t,z_i+c_i(t)-y(t)\right)_0\ge 0,\ \sum_{i=1}^k  \left(\dot y_i^t,z_i+c_i(t)-x(t)\right)_0\ge 0,\    z_i\in C_i,
\end{equation} so that the arguments of Step~1 apply 
to
$$
  f(h)=\sum_{i=1}^k\left(\dot x_i^t,y(t+h)-y(t)+c_i(t)-c_i(t+h)\right)_0
$$
(similarly for the second sum of (\ref{weclaim}) with $z_i=x(t+h)-c_i(t+h)$)
to show that the sums in (\ref{vanish2}) vanish. To establish (\ref{weclaim}), we first rewrite it as
$$
\sum_{i=1}^k\left(\dot x_i^t,z_i+c_i(t)\right)_0-\left(\dot x(t),y(t)\right)_0\ge 0,\ \ \sum_{i=1}^k  \left(\dot y_i^t,z_i+c_i(t)\right)_0-\left(\dot y(t),x(t)\right)_0\ge 0,
$$
and then 
prove that
\begin{equation}\label{establ}
\left(\dot x(t),y(t)\right)_0=\left(\dot x(t),x(t)\right)_0\ \ {\rm and} \ \ \left(\dot y(t),x(t)\right)_0=\left(\dot y(t),y(t)\right)_0,\ \ t\in[0,T],
\end{equation}
so that (\ref{weclaim}) becomes a consequence of (\ref{dnc}).
To prove (\ref{establ}) we use \eqref{eq:constantDistance} and observe that
$$
0=\frac{d}{dt}\|x(t)-y(t)\|^2_0=
-\left(\dot x(t),y(t)-x(t)\right)_0-\left(\dot y(t),x(t)-y(t)\right)_0
$$
But $x(t),y(t)\in C(t)$ and both these functions are solutions of sweeping process (\ref{swgen}). Therefore,  
$\left(\dot x(t),y(t)-x(t)\right)_0\geqslant 0$ and $\left(\dot y(t),x(t)-y(t)\right)_0\geqslant 0,$
which implies (\ref{establ}).

\vskip0.2cm

\noindent The proof of the lemma is complete. \qed

\vskip0.2cm \noindent We acknowledge that the idea of the proof of Step 1 of Lemma~\ref{th:constDistance} has been earlier used by Krejci in the proof of \cite[Theorem 3.14]{Krejci1996}, which would suffice for the proof when $k=1.$ The achievement of Lemma~\ref{th:constDistance} is in considering $k>1$, thus the new Step~2. Accordingly, the proof of the next theorem follows the lines  of \cite[Theorem 3.14]{Krejci1996} with Lemma~\ref{th:constDistance} used to justify (\ref{onlynew}),
   which is the place of the proof that needed further arguments when moving to $k>1$. We present a proof for completeness (Appendix~\ref{appen}) also because \cite{Krejci1996}
employs slightly different notations. The theorem effectively states that any bounded solution of a $T$-periodic sweeping process is asymptotically $T$-periodic, which facts is known in differential equations as Massera's theorem \cite{massera}.  

\vskip0.2cm

\begin{theorem} \textit{\textbf{ {\bf (}Massera-Krejci theorem for sweeping processes with a moving set of the form \boldmath{$C(t)=\cap_{i=1}^k (C_i+c_i(t))$}{\bf )} }}
\label{lemattract}
Let $t\mapsto C(t)$ be a set-valued uniformly bounded function of the form  (\ref{eq:finiteIntersectionC})-(\ref{nonempty}) with convex closed  $C_i$ and Lipschitz-continuous single-valued $T$-periodic $c_i(t)$. Then the set $X(t)$ of $T$-periodic solutions of (\ref{swgen}) is a global attractor of (\ref{swgen}). 
\end{theorem}

\noindent  The proof of Theorem~\ref{lemattract} is given in Appendix~\ref{appen}.

\vskip0.2cm

\noindent An interested reader can note that sweeping process (\ref{swgen}) with $k=1$ converts to a perturbed sweeping process 
$
  -\dot \xi=N_{C_1}^0(\xi)+\dot c_1(t)
$
with an immovable constraint by the change of the variables $\xi(t)=y(t)-c_1(t)$,    while it is not clear whether or not (\ref{swgen})  converts to a perturbed sweeping process with a constant constraint when $k>1.$ This further highlights the difference between the cases $k=1$ and $k>1$ as long as potential alternative methods of analysis of the dynamics of (\ref{swgen}) are concerned.

\subsection{Strengthening of the conclusion of section \ref{sec41} in the case of a moving constraint given by a polyhedron with translationally moving facets}\label{strengthening}

\noindent When applied to a one-dimensional network of elastoplastic springs (\ref{eq1})-(\ref{eq5}), the existence of a periodic attractor $X(t)$ for the associated sweeping process (\ref{sw1}) follows from Theorem~\ref{thm1}. A new geometric property of $X(t)$ that comes with considering the sweeping process (\ref{sw1}) is due to the 
polyhedral shape of the moving constraint  $\Pi(t)\cap V$, see Section~\ref{polyhedronsection}. Theorem~\ref{thmnew} below states that even if $X(t)$ consists of several periodic solutions, they all exhibit certain identical behavior. 

\vskip0.2cm

\noindent As earlier, let $E$ be a finite-dimensional linear vector space equipped with a scalar product $(\cdot,\cdot)_0$ and let ${\rm ri}(X)$ be the relative interior of the convex set $X\in C([0,T],E)$.

\begin{theorem}\label{thmnew} Assume that a uniformly bounded set-valued function $t\mapsto C(t)$ is given by 
\begin{equation}\label{simp0}
     C(t)=\bigcap\limits_{i=1}^m\left\{x\in E:c_i^-(t)\le \left(n_i,x\right)_0\le  c_i^+(t)\right\},\quad t\geqslant 0,
\end{equation}
where $c_i^-,c_i^+$ are single-valued globally Lipschitz continuous functions, $n_i$ are given vectors from $E$. Then the set $X(t)$ of $T$-periodic solutions of sweeping process (\ref{swgen}) is the global attractor of (\ref{swgen}). Furthermore, 
$X\subset C([0,T],E)$ is closed and convex, and all the interior solutions of $X$ follow the same pattern of motion in the sense that 
\begin{equation}\label{Jprop}
 J(t,x(t))= J(t,y(t)),\ { for\ all\ }x,y\in {\rm ri}(X),\  t\ge 0,
\end{equation}
where  $J(t,x)$ is the active set of the polyhedron $C(t)$ given by 
$$
J(t,x)=\left\{i\in\overline{-m,-1}:(n_{-i},x)_0=c_{-i}^-(t)\right\}
\cup\left\{i\in\overline{1,m}:(n_i,x)_0=c_i^+(t)\right\}.
$$
\end{theorem}

\vskip0.2cm

\noindent Theorem~\ref{thmnew} is a corollary of Theorem~\ref{thm1} except for the property (\ref{Jprop}) which comes from the polyhedral shape of the moving constraint $C(t).$ The property (\ref{Jprop}) follows from the following general result.

\begin{lemma} \label{lemma:sameFacesPerSolutions} 
 Consider an arbitrary convex set $B$  embedded into a convex polyhedron:
$$
   B\subset\bigcap_{i=1}^{k}\left\{x\in E:\left( y_i,x\right)_0\le  b_i\right\},
$$
where $y_i\in E$ and $ b_i\in\mathbb{R}.$ If $x_1,x_2\in {\rm ri}(B),$ then, for all $i\in\overline{1,k},$
\begin{equation}\label{implication}
   \left( y_i,x_1\right)_0=  b_i\quad\mbox{if and only if}\quad  \left( y_i,x_2\right)_0= b_i.
\end{equation}
\end{lemma}
\noindent{\bf Proof.} Consider 
$$
  x_\theta=\theta x_1+(1-\theta)x_2.
$$
Since $x_1,x_2\in {\rm ri}(B)$, there exists $\varepsilon>0$ such that $x_{-\varepsilon}\in {\rm ri} (B)$ and $x_{1+\varepsilon}\in{\rm ri} (X).$ Put $\bar x_1=x_{-\varepsilon},$ $\bar x_2=x_{1+\varepsilon}.$ Then there exist 
$\theta_1,\theta_2\in(0,1),\, \theta_1\neq \
\theta_2$, such that 
\begin{equation}\label{twoformula}
x_1=\theta_1 \bar x_1+(1-\theta_1)\bar x_2,\quad
x_2=\theta_2 \bar x_1+(1-\theta_2)\bar x_2.
\end{equation}

\noindent Assume that 
$\left( y_i, x_1\right)_0=b_i$. Then
\begin{equation}\label{canhold}
\theta_1 (\left(  y_i,\bar x_1\right)_0- b_i)=-(1-\theta_1)(\left(  y_i,\bar x_2\right)_0- b_i)
\end{equation}
by using the first formula of (\ref{twoformula}). Since $\bar x_1,\bar x_2\in B$, one has $\left(y_i,\bar x_1\right)_0- b_i\le 0$ and $\left(y_i,\bar x_2\right)_0- b_i\le 0$. Therefore formula (\ref{canhold}) can only hold  when both 
$\left( y_i,\bar x_1\right)_0- b_i$ and $\left( y_i,\bar x_2\right)_0-b_i$ vanish.
Hence 	
$$\left( y_i,x_2\right)_0=\theta_2\left( y_i,\bar x_1\right)_0+(1-\theta_2)\left( y_i,\bar x_2\right)_0=\theta_2 b_i+(1-\theta_2) b_i= b_i.$$
The reverse implication in (\ref{implication})  can be proved by analogy.\qed

\subsection{Application: an analytic condition for the convergence of the stresses of elastoplastic systems to an attractor}\label{sec1a}

\noindent Let $J_C(x)$ be the active set of the parallelepiped $C$, i.e. 
$$
   J_C(x)=\left\{i\in\overline{-m,-1}:x_{-i}=c_{-i}^-\right\}\cup \left\{i\in\overline{1,m}:x_i=c_i^+\right\}.
$$
 A direct consequence of Theorem~\ref{thmnew} is the following result about asymptotic behavior of the stresses of  the elastoplastic system (\ref{eq1})-(\ref{eq5}).

\begin{theorem}\label{cor0} 
Let the conditions of Theorem~\ref{moreauthm} hold and both displacement-controlled and stress-controlled loadings are $T$-periodic. 
Then, for any initial condition at $t=0,$ the stresses $s_1(t),...,s_m(t)$ of the springs 
converge, as $t\to\infty,$ to the attractor
$$
  S(t)=A\left(X(t)-h(t)+g(t)\right),
$$
where $X(t)$ is the set of all $T$-periodic solutions of sweeping process (\ref{sw1}), and $h(t)$ and $g(t)$ are the effective loadings given by \eqref{h} and \eqref{g(t)}.
The functions of $S(t)$ have equal derivatives for a.a. $t\ge 0$ as per (\ref{equivV}) and, moreover,
\begin{equation}
 J_C\left(\bar s_i(t)\right)= J_C\left(\hat s_i(t)\right),\ { for\ all\ }\bar s,\hat s\in {\rm ri}(S),\  t\ge 0.
\label{eq:sameIndices}
\end{equation} 
\end{theorem}
 
\noindent {\bf Proof.} We apply Theorem~\ref{thmnew} with $C(t)=\Pi(t)\cap V,$ where $\Pi(t)$ and $V$ are those defined in Theorem~\ref{moreauthm}. Since $\Pi(t)$ is uniformly bounded in $t\in[0,T]$, same holds for $C(t).$ Thus, the conditions of  
Theorem~\ref{thmnew} are satisfied with $c_i^-(t)=c_i^-+a_i h_i(t)$ and $c_i^+(t)=c_i^++a_i h_i(t),$ and   Theorem~\ref{thmnew} implies that
$$
  J(t,A^{-1}\bar s(t)+h(t)-g(t))=J(t,A^{-1}\hat s(t)+h(t)-g(t)), 
\ {\rm for\ all\ }\bar s,\hat s\in {\rm ri}(S),\  t\ge 0,
$$
which equivalent formulation is (\ref{eq:sameIndices}).  Other statements of Theorem~\ref{cor0} follow from Theorem~\ref{thmnew} just directly. The proof of the theorem is complete.\qed

\vskip0.2cm

\noindent Property (\ref{eq:sameIndices}) says that, for any $i\in \overline{1,m}$, the spring $i$ will asymptotically execute a certain pattern of elastoplastic  deformation which doesn't depend on the state of the network at the initial time. 

\vskip0.2cm

\noindent We remind the reader that if $s_i,$ $i\in\overline{1,m},$ are the stresses of springs, then the  quantities $\dfrac{1}{a_i}s_i(t)$ that appear in (\ref{eq:sameIndices}) are the elastic elongations of the springs.

\vskip0.2cm

\noindent {\bf Example (continued).} For the elastoplastic system (\ref{eq1})-(\ref{eq5}) of Fig.~\ref{ex1fig} with $T$-periodic displacement-controlled and stress-controlled loadings $l(t)$ and $h(t)$, Theorem~\ref{cor0} implies the convergence of stresses $s(t)$ to a $T$-periodic attractor $S(t)$ provided that property (\ref{ex1_safeload}) holds. Furthermore, the functions of $A^{-1}S(t)+h(t)-g(t)$ are all non-constant, if  (\ref{ex1_nonconstant}) is satisfied. In the next section of the paper we offer a general result which will, in particular, imply that the attractor $A^{-1}S(t)+h(t)-g(t)$ consists of a single solution. 

\section{Stabilization to a unique non-stationary periodic solution}\label{sec5}

In this section we first prove that the periodic attractor $X(t)$ of a general sweeping process (\ref{swgen}) in a vector space of dimension $d$ consists of just one non-stationary $T$-periodic solution, when  
the normals of any $d$ different facets of the moving polyhedron $C(t)$ are linearly independent.  Then we give a sufficient condition for such a requirement to hold for the sweeping process (\ref{sw1}) coming from the elastoplastic system (\ref{eq1})-(\ref{eq5}).

\subsection{Stabilization of a general sweeping process with a polyhedral moving set}\label{sec51}

As earlier, let $E$ be a linear vector space of dimension $d$ and let $(\cdot,\cdot)_0$ be a scalar product in $E.$
\vskip0.2cm

\noindent
In this subsection it will be convenient to rewrite the set (\ref{simp0}) in the following form 
\begin{equation}\label{simp1}
     C(t)=\bigcap\limits_{i=1}^k\left\{x\in E: (x,n_i)_0\le c_i(t)\right\},\quad t\geqslant 0,
\end{equation}
where $c_i$ are single-valued functions and $n_i$ are given vectors of $E$. 
 The advantage of form (\ref{simp1}) compared to (\ref{simp0}) is that any vector of $N_{C(t)}^0(x)$ has non-negative coordinates in the basis formed by the normals $n_1,...,n_k,$ as our Lemma~\ref{lemma:conicDecomposition} shows. Then we establish the following result about global asymptotic stability of sweeping processes.
 
\begin{theorem}\label{thm1.5} Let  $t\mapsto C(t)$ be a uniformly bounded set-valued function given by \eqref{simp1}, where the functions $c_i$ are globally Lipschitz continuous and $k\ge \dim E.$ Assume that any $\dim E$ vectors out of 
the collection $\{n_i\}_{i=1}^k\subset E$ are linearly independent and the cardinality of the set 
$$J(t,x)=\{i\in\overline{1,k}:(x,n_i)_0=c_i(t)\}$$
 doesn't exceed $\dim E$ for all $x\in C(t)$ and $t\in[0,T]$. Then
$X(t)$ contains at most one non-constant $T$-periodic solution.
\end{theorem}

\noindent Note, $n_i$ in (\ref{simp1}) are, generally speaking, different from $n_i$ in (\ref{simp0}), but we use same notation as it shouldn't cause confusion. Accordingly, the active set $J(t,x)$ of Theorem~\ref{thm1.5} is different from the active set $J(t,x)$ of Theorem~\ref{thmnew}.

\begin{lemma}
Assume, that for each $t\in[0,T]$ and $x\in C(t)$ the collection of vectors $\{n_i:i\in J(t,x)\}$ is linearly indepentent. Then for a solution $x(t)$ of sweeping process (\ref{swgen}) there is a collection of integrable non-negative $\lambda_i:[0,T]\to \mathbb{R}^+ \cup\{0\}$, $i\in\overline{1,k},$  such that 
\begin{equation}
-\dot x(t)=\sum_{i=1}^k\lambda_i(t)\,n_i,\quad  \mbox{for a.a. } t\in [0,T]. 
\label{eq:betaDecomposition}
\end{equation}
\label{lemma:conicDecomposition}
\end{lemma}

\noindent {\bf Proof.} 
Recall, that for a fixed $t\in[0,T],$ the normal cone 
(\ref{normal_cone_standard})
to the set   $C(t)$ of form (\ref{simp1}) can be equivalently formulated as  
\begin{equation}
 \renewcommand{\arraystretch}{1.8}
N^0_{C(t)}(x)=\left\{\begin{array}{ll}
\{0\}, & \  \ {\rm if \ } x\in {\rm int}\,C(t),\\
\left\{\left(\sum\limits_{i\in J(t,x)} \lambda_i\,n_i\right): \lambda_i\ge 0\right\}, & \  \ {\rm if\ } x\in\partial C(t),\\
\emptyset, & \ \ {\rm if}\ x\not\in C(t).
\end{array}
\right.
\label{eq:normlConePolyhedral}
\end{equation}
Here ${\rm int}\, C(t)$ and $\partial C(t)$ are respectively the interior and the boundary of  $C(t).$
Therefore, for a.a. fixed $t\in[0,T],$ the existence of $\lambda(t)\in \mathbb{R}^k$, $\lambda_i(t)\ge 0,$ $i\in\overline{1,k},$ verifying
\begin{equation}
-\dot x(t)=\sum_{i\in J(t,x(t))}\lambda_i(t) n_i
\label{eq:dotIsSum}
\end{equation}
follows from the inclusion (\ref{swgen}).  We set $\lambda_i(t)=0,$ if $i\not\in J(t,x(t)).$ The proof of Lebesgue measurability of $\lambda(t)$ will be split into several steps.

\vskip0.2cm

\noindent {\bf Step 1.} First we observe that, for any $\hat t\in[0,T]$,
$$
  \mbox{the set }T_{\hat t}=\left\{t\in[0,T]:J(t,x(t))=J(\hat t,x(\hat t))\right\}\mbox{ is measurable}.
$$
This follows from the fact that the set 
$
  \left\{t\in[0,T]:\left<x(t),n_i(t)\right>-c_i(t)=0\right\}$ is measurable  for each fixed index $i\in\overline{1,k}$ and that $J(\hat t,x(\hat t))\subset \overline{1,k}.$

\vskip0.2cm

\noindent {\bf Step 2.} Now we fix some $\hat t\in[0,T]$  and prove that, for any Borel set $B\subset\mathbb{R}^k,$ 
$$
  \mbox{the set }T_{\hat t}(B)=\left\{t:[0,T]:\lambda(t)\in B,\ J(t,x(t))=J(\hat t,x(\hat t))\right\} \mbox{ is measurable.}
$$ 
If inclusion (\ref{swgen}) doesn't hold at $\hat t$ and ${\rm mes}(T_{\hat t})=0,$ then $T_{\hat t}(B)$ is measurable and ${\rm mes}(T_{\hat t}(B))=0.$ If (\ref{swgen}) doesn't hold  at  $\hat t$ and ${\rm mes}(T_{\hat t})>0,$ then we can find $\tilde t\in T_{\hat t}$ such that (\ref{swgen}) does hold at $\tilde t.$ Since $T_{\hat t}=T_{\tilde t}$, we conclude that one won't  restrict generality of the proof, if assume that (\ref{swgen}) holds for the initially chosen $\hat t\in[0,T].$ 

\vskip0.2cm

\noindent  Let $\bar n_1,...,\bar n_d$ be any basis in $E$ such that
$$
   \bar n_i=n_i,\quad\mbox{for all}\ i\in J(\hat t,x(\hat t)),
$$ therefore it depends on $\hat t.$
Denote by  $S_{\hat t}:E\to\mathbb{R}^d$ the bounded linear map which maps every vector from $E$ to its coordinates in terms of $\{\bar n_i\}_{i=1}^d$. Then \eqref{eq:dotIsSum} necessarily means that  
$$
   \lambda(t)=-S_{\hat t}\dot x(t),\quad\mbox{for a.a. }t\in[0,T]\mbox{ such that }J(t,x(t))=J(\hat t,x(\hat t)).
$$
Therefore, up to a subset 	of $[0,T]$ of zero measure, 
$$ T_{\hat t}(B) =\left\{t\in[0,T]:-S_{\hat t}\dot x(t)\in B,\ J(t,x(t))=J(\hat t,x(\hat t))\right\}= (-S_{\hat t}\dot x)^{-1}(B)\cap T_{\hat t},
$$
and the measurability of $T_{\hat t}(B)$ follows by combining the continuity of $x$ and the conclusion of Step~1.

\vskip0.2cm

\noindent {\bf Step 3.} We finally fix a Borel set $B\subset\mathbb{R}^k$ and prove the measurability of the set
\begin{equation}\label{La}
   \lambda^{-1}(B)=\left\{t\in[0,T]:\lambda(t)\in B\right\}.\end{equation}
Since $J(\hat t,x(\hat t))$ can take only a finite number of (set-valued) values when $\hat t$ varies from $0$ to $T,$ then there is a finite sequence $t_1,...,t_K\in[0,T]$ such that
$$
   [0,T]=\bigcup_{i\in\overline{1,K}}T_{t_i},
$$
and so we can rewrite (\ref{La}) as follows
$$
  \lambda^{-1}(B)=\bigcup_{i\in\overline{1,K}}T_{t_i}(B),
$$
which is a finite union of measurable sets. The proof of the measurability of  $\lambda$ is complete.

\vskip0.2cm

\noindent The integrability of $\lambda$ on $[0,T]$ now follows from its boundedness. Indeed, since,  $\|\dot x(t)\|_0\le M$ for a.a. $t\in[0,T]$ and some $M>0$ \cite[p.13]{kunze}, one has 
$$|\lambda_i(t)|\le\|\lambda(t)\|=M\max_{i\in\overline{1,K}}\|S_{t_i}\|,\quad {\rm for \ a.a.\ }t\in[0,T].$$
The proof of the lemma is complete.\qed

\vskip0.2cm

\noindent {\bf Proof of Theorem~\ref{thm1.5}.} Let $x(t)$ and $y(t)$ be two non-constant distinct $T$-periodic solutions of (\ref{sw1}). Theorem~\ref{thmnew} implies that we won't lose generality by assuming that 
\begin{equation}\label{JJ}
   J(t,x(t))=J(t,y(t)).
\end{equation}
When applying Theorem~\ref{thmnew}
we used the fact that the set (\ref{simp1}) can be expressed in the form (\ref{simp0}) due to  the uniform boundedness of  $C(t)$.

\vskip0.2cm

\noindent  The proof is by reaching a contradiction with the fact that $x(t)$ and $y(t)$ are distinct.

\vskip0.2cm

\noindent By replacing $-\dot x(t)$  
by its representation from Lemma \ref{lemma:conicDecomposition}, one gets 
\begin{equation}\label{prope}
0=\int_0^\top-\dot x(t)dt=
\int_0^\top \sum_{i=1}^k \lambda_i(t)n_i dt=
\sum_{i=1}^k \int_0^\top \lambda_i(t)dt\, n_i,
\end{equation}
where $\lambda_i(t)\ge 0.$ Since $x(t)$ is non-constant, the set  
$$\hat J:=\left\{i\in\overline{1,k}:\int_0^\top \lambda_i(t)dt> 0\right\}$$
is non-empty. The following two cases can take place.

\vskip0.1cm

\noindent {\bf 1)} $\{n_i:i\in \hat J\}$ is a linearly independent system. But property (\ref{prope}) yields 
$$\sum_{i\in \hat J} \int_0^\top \lambda_i(t)dt\, n_i=0,
$$
that, for linearly independent vectors $n_i$, can happen only when $\int_0^\top \lambda_i(t)dt\equiv0,$ $i\in\overline{1,k}.$ Therefore case 1) cannot take place as $x(t)$ is non-constant.

\vskip0.1cm

\noindent {\bf 2)} The vectors of $\{n_i:i\in \hat J\}$ are linearly dependent. Since, by the assumption of the theorem, any $d$ vectors from $\{n_i:i\in \hat J\}$ are linearly independent, one must have $|\hat J|>d$. Let us show this leads to a contradiction as well.

\vskip0.1cm

\noindent Since for each $j\in \hat J,$ the function $ \lambda_j(t)$ is positive on a set of positive measure, there are time moments $\{t_j\}_{j\in\hat J}$, where (\ref{JJ}) holds along with
$$-\dot x(t_j)=\sum\limits_{i=1}^ k  \lambda_i(t_j)n_i \quad{\rm and}\quad \lambda_j(t_j)>0. $$
This and \eqref{eq:normlConePolyhedral} imply 
$$
   j\in J(t_j,x(t_j))\quad \text{ and  by \eqref{JJ}} \quad j\in J(t_j,y(t_j)),\quad j\in\hat J,
$$
or, equivalently,
$$
\left(x(t_j),n_j\right)_0=c_j(t_j)\quad{\rm and}\quad \left(y(t_j),n_j\right)_0=c_j(t_j),\quad j\in \hat J.
$$
Therefore,
$$
\left(x(t_j)-y(t_j),n_j\right)_0=0,\quad j\in\hat J,
$$
and, by Lemma~\ref{th:constDistance},
$$
\left(x(0)-y(0),n_j\right)_0=0,\quad j\in\hat J.
\label{eq:niEquality}
$$
But $|\hat J|>d$  and so $\{n_i:i\in \hat J\}$ contains $d$ linearly independent vectors, which form a basis of $\mathbb{R}^d.$ Therefore, $x(0)=y(0)$, which is a contradiction.
\qed

\vskip0.2cm

\noindent Theorem~\ref{thm1.5} can be used for stabilization of general sweeping process with polyhedral moving set such as those considered e.g. in Colombo et al \cite{colombo} and Krejci-Vladimirov \cite{krejci1}. 

\vskip0.2cm

\noindent A fundamental case where Theorem~\ref{thm1.5} allows to stabilize an elastoplastic system (\ref{eq1})-(\ref{eq5}) to a single periodic solution is when $V$ cut $\Pi(t)$ along a simplex. Testing the set $\Pi(t)\cap V$ for being a simplex can be executed for any given  elastoplastic system (\ref{eq1})-(\ref{eq5}) using the algorithms of computational geometry (e.g. Bremner et al \cite{bremner} can be used to compute the vertexes of $\Pi(t)\cap V$ whose number needs to equal $m+1$). 

\vskip0.2cm

\noindent At the same time, establishing analytic criteria for stabilization to occur could be of great use in materials science. A simple criterion of this type is offered in the next section of the paper.

\subsection{Application: an analytic condition for stabilization of elastoplastic systems to a unique periodic regime}\label{sec52}

\noindent Next theorem is the main result of this paper. It can be viewed as an analogue of high gain feedback stabilization in control theory. Indeed, one of the two central assumptions of the theorem is $q=n-2$, which means that the elastoplastic system has a sufficient number of control variables to be fully controllable and thus stabilizable. The second central assumption is assuming that the magnitude of the stress-controlled loading  is high enough which literally resembles the high gain requirement of feedback control theory.

\vskip0.2cm

\noindent The idea of Theorem~\label{mainthm} is based on a simple fact that the moving parallelepiped $\Pi(t)$ intersects the plane $V$ along a simplex, if the the plane $V$ is close to the vertex of the parallelepiped, see Fig.~\ref{ex1fig_}(d). At the same time, this geometric statement turned out to hold only if $q=n-2$.

\begin{theorem}\label{mainthm} In the settings of Proposition~\ref{cor2} assume that the stress loading $h(t)$ is large in the sense that 
\begin{equation}\label{bigforce}
\left<\bar u,\left(\begin{array}{c}
\bar c^{k}_1\\ \vdots\\ \bar c^k_{j-1} \\ \bar c^{-k}_j \\ \bar c^k_{j+1}\\ \vdots \\
\bar c^k_m
\end{array}\right)+Ah(t)
\right>\cdot \left<\bar u, \bar c^k+Ah(t)\right>\le 0,
\quad j\in\overline{1,m},\ t\in[0,T],
\end{equation}
holds for at least one $k\in\{-1,+1\}$. Further assume that the displacement-controlled loading $g(t)$ is large in the sense of (\ref{formulaprop3}). Then, there exists a $T$-periodic function $s_0(t)$ such that $\|s(t)-s_0(t)\|\to 0$ as $t\to\infty$ for the stress component $s(t)$ of any solution of the quasistatic evolution problem (\ref{eq1})-(\ref{eq5}).
\end{theorem}

\begin{remark}\label{rem5} Following the lines of Remark~\ref{remark4}, we consider the left-hand-side of (\ref{bigforce}) as a polynomial $P\left(\left<\bar u,Ah(t)\right>\right)$ in $\left<\bar u,Ah(t)\right>$, so that the branches of the polynomial are pointing upwards. Therefore, condition (\ref{bigforce}) is the requirement for $\left<\bar u,Ah(t)\right>$ to stay between the roots of the polynomial. Note, one root of $P\left(\left<\bar u,Ah(t)\right>\right)$ is given by $\left<\bar u,Ah(t)\right>=-\left<\bar u,\bar c^k\right>$. By computing the derivative $P'\left(-\left<\bar u,\bar c^k\right>\right)$  one concludes that $\left<\bar u,Ah(t)\right>=-\left<\bar u,\bar c^k\right>$ is the smaller or larger root of $P\left(\left<\bar u,Ah(t)\right>\right)$ according to whether $k=+1$ or $k=-1.$ Therefore, a sufficient condition for (\ref{bigforce}) to hold with $k=+1$ and $k=-1$ are
$$
\left<\bar u,-\bar c^+\right>\le
\left<\bar u,Ah(t)\right>\le\min_{j\in\overline{1,m}}
\left<\bar u,-\left(\begin{array}{c}
\bar c^{+}_1\\ \vdots\\ \bar c^+_{j-1} \\ \bar c^{-}_j \\ \bar c^+_{j+1}\\ \vdots \\
\bar c^+_m
\end{array}\right)
\right> 
$$ 
and
$$
\max_{j\in\overline{1,m}}
\left<\bar u,-\left(\begin{array}{c}
\bar c^{-}_1\\ \vdots\\ \bar c^-_{j-1} \\ \bar c^{+}_j \\ \bar c^-_{j+1}\\ \vdots \\
\bar c^-_m
\end{array}\right)
\right>\le \left<\bar u,Ah(t)\right>
\le \left<\bar u,-\bar c^-\right> 
$$ 
respectively.
\end{remark}

\noindent {\bf Proof.} We are going to prove that the conditions of Theorem~\ref{thm1.5} hold for the sweeping process \eqref{eq:sw3} of the elastoplastic system given. Since the set 
\begin{equation}\label{333}
  C(t)=
\bigcap_{i=1}^{m}\left\{x\in V:c_i^-\hskip-0.05cm+\hskip-0.05cma_ih_i(t)\hskip-0.05cm\le\hskip-0.05cm a_i x_i\hskip-0.05cm\le\hskip-0.05cm c_i^+\hskip-0.05cm+\hskip-0.05cma_ih_i(t)\right\}
\end{equation}
is just a parallel displacement of the polyhedron (\ref{222}) of sweeping process \eqref{eq:sw3}, it is sufficient to prove that conditions of Theorem~\ref{thm1.5} hold for the set (\ref{333}). More precisely, we prove that 
conditions of Theorem~\ref{thm1.5} hold for the set (\ref{333}) after it is expressed in the form (\ref{simp1}).

\vskip0.2cm
 
\noindent Fix $t\in[0,T]$ and $j\in\overline{1,m}.$ Denote by $\xi^j$ the solution of the system of $m$ equations
\begin{eqnarray}
&&       \left<\bar u,A\xi^j\right>=0,\label{xi1}\\
    &&   a_i\xi_i^j=\bar c_i^k+a_i h_i(t),\quad i\in\overline{1,m}, i\neq j. \label{xi2}
\end{eqnarray}
The solution $\xi^j$ is unique by Lemma~\ref{lemu0}  used for $x=\bar u$. 
Observe, that 
\begin{multline}
\la\bar u, A\left(A^{-1}\bar c^k+h(t)\right) \ra=0\quad \Leftrightarrow \quad \xi^j=\bar c^k,j\in\overline{1,m}\quad \Leftrightarrow \\ \Leftrightarrow \quad C(t)=\{A^{-1}\bar c^k+h(t)\}.
\label{eq:singletonCcondition}
\end{multline}
Indeed, assume that there exists $x\in C(t)$ such that $x\not=A^{-1}\bar c^k+h(t).$ Then $x$ can be expressed as $x=A^{-1}c+h(t)$ for some $c\in C.$ On the other hand, $x\in V$ implies $\left<\bar u,Ax\right>=0.$ Therefore, $\left<\bar u,\bar c^k- c\right>=0$ and by just expanding the scalar product we get the existence of two 
indices $i_1,i_2\in \overline{1,m}$, such that
$$
\bar u_{i_1}(\bar c_{i_1}^k-c_{i_1})>0, \qquad \bar u_{i_2}(\bar c_{i_2}^k-c_{i_2})<0,
$$
which is impossible by the construction of $\bar c^k$. Therefore $C(t)$ is a singleton, if (\ref{eq:singletonCcondition}) holds. But if $C(t)$ is a singleton for at least one $t\in[0,T]$, then the statement of the theorem becomes trivial. That is why we now focus on the case where \eqref{eq:singletonCcondition} doesn't hold on $[0,T].$ Below we will complete the proof of \eqref{eq:singletonCcondition} by showing that all $\xi^j\in C(t)$.

\vskip0.2cm

\noindent In what follows, we show that the conditions of Theorem~\ref{thm1.5} hold by proving that the points $\xi^j$, $j\in\overline{1,m},$ are vertices of a $m-1$-simplex that coincides with $C(t)\subset V$(recall that ${\rm dim}\, V=m-{\rm dim}\, U=m-1$ in this case).

\vskip0.2cm
\noindent {\bf Step 1:} {\it It holds $\xi^j\in C(t)$, $j\in\overline{1,m}.$} Based on formula (\ref{333}), we  have to show that 
\begin{equation}\label{444}
  c_j^-+a_jh_j(t)\le a_j\xi_j^j\le c_j^++a_jh_j(t),\quad j\in\overline{1,m}.
\end{equation}
Fix $j\in\overline{1,m}$ and consider the function
$$
  b(x)=
\left<\bar u,\left(\begin{array}{c}
\bar c^{k}_1\\ \vdots\\ \bar c^k_{j-1} \\ x \\ \bar c^k_{j+1}\\ \vdots \\
\bar c^k_m
\end{array}\right)+\left(\begin{array}{c}
a_1h_1(t)\\ \vdots\\ a_{j-1}h_{j-1}(t) \\ 0 \\  a_{j+1}h_{j+1}(t)\\ \vdots \\
a_m h_m(t)
\end{array}\right)
\right>.
$$
By the definition, $a_j\xi_j^j$ is the unique root of the equation $b(x)=0$. On the other hand, condition (\ref{bigforce}) implies that $b(\bar c_j^{-k}+a_jh_j(t))\cdot b(\bar c_j^k+a_jk_j(t))\le0$, so that the unique zero of $b(x)$ must be located between the numbers $\bar c_j^{-k}+a_jh_j(t)$ and $\bar c_j^k+a_jk_j(t)$.

\vskip0.2cm

\noindent {\bf Step 2:} {\it The vertices $\xi^j$, $j\in\overline{1,m},$ form an $m-1$-simplex.} For a given $j\in\overline{1,m},$ we need to show that $m-1$ vectors 
$$
  \zeta^{ij}=\xi^i-\xi^j,\quad i\in\overline{1,m},\ i\not=j,
$$
are linearly independent. From (\ref{xi1}) we have 
$$
   \left<\bar u,A\zeta^{ij}\right>=0
$$
while from (\ref{xi2}) we get 
\begin{equation}
  \zeta^{ij}=(0,...,0,\zeta_i^{ij},0,...,0,\zeta_j^{ij},0,...,0)^\top.
\label{eq:zetaForm}
\end{equation}
Combining these two properties we conclude that
$$
u_i a_i\zeta_i^{ij}+u_ja_j\zeta_j^{ij}=0.
$$
By Lemma~\ref{lemu0} $u_ia_i\neq0$ and $u_j a_j\neq 0$, therefore we either have $\zeta_i^{ij}=\zeta_j^{ij}=0$ or 
$\zeta_i^{ij}\zeta_j^{ij}\not=0.$ Observe that the former case is impossible. 
Indeed, if $\xi^{j_1}=\xi^{j_2}$ for some $j_1\not=j_2,$ then (\ref{xi2}) implies $\xi^{j_1}=\xi^{j_2}=A^{-1}(\bar c^k+Ah(t))$, which leads to (\ref{eq:singletonCcondition})  when plugged to (\ref{xi1}) which we already excluded.

\vskip0.2cm

\noindent It remains to notice that property $\zeta_i^{ij}\zeta_j^{ij}\not=0,$ $i\not=j$ implies that the vectors (\ref{eq:zetaForm}) are linearly independent through $i\not=j,$ $i\in\overline{1,m}.$

\vskip0.2cm

\noindent {\bf Step 3:} {\it We claim that $C(t)={\rm conv}\left\{\xi^j, \ j\in\overline{1,m}\right\}.$} From {Step~1}, $C(t)\supset{\rm conv}\left\{\xi^j, \ j\in\overline{1,m}\right\},$ so it remains to prove that $C(t)\subset {\rm conv}\left\{\xi^j, \ j\in\overline{1,m}\right\}.$

\vskip0.2cm

\noindent We fix $\hat j\in \overline{1,m}$ and consider a facet ${\rm conv}\left\{\xi^i,\ i\not=\hat j,\right\}$  of the simplex ${\rm conv}\left\{\xi^i,\  i\in\overline{1,m}\right\}$. Observe from \eqref{xi2} that all  vertices of the facet share their $\hat j$-th coordinate. Therefore the whole facet belongs to the plane
$$
  L^{\hat j}=\{x\in V:a_jx_j=\bar c_{\hat j}^k+a_{\hat j} h_{\hat j}(t)\}.
$$
Therefore,
$$
{\rm conv}\left\{\xi^i,\  i\in\overline{1,m}\right\}=\bigcap_{j=1}^m \left\{x\in V: p_j\left(a_jx_j-\bar c_j^k-a_j h_j(t)\right)\le 0\right\},
$$
where $p_j\in\{-1,1\}$ are suitable signs. On the other hand, by (\ref{333}),
\begin{equation}\label{takes}
  C(t)\subset \bigcap_{j=1}^m \left\{x\in V: q_j\left(a_jx_j-\bar c_j^k-a_j h_j(t)\right)\le 0\right\},
\end{equation}
where $q_j\in\{-1,1\}$ are suitable signs. Since by Step~2, ${\rm conv}\{\xi_j,\ j\in\overline{1,m}\}\subset C(t)$, we get $p_j=q_j,$ $j\in\overline{1,m}.$ But then (\ref{takes}) takes the form $C(t)\subset 
{\rm conv}\{\xi_j,\ j\in\overline{1,m}\}.$

\vskip0.2cm

\noindent The proof of the theorem is complete. \qed

\vskip0.2cm

\noindent{\bf Example (continued).} Applying Theorem~\ref{mainthm} to the elastoplastic system of Fig.~\ref{ex1fig} (where we have $q=n-2$) we use earlier formulas (\ref{UbasisEx1}) and (\ref{ex1form}) together with Remark~\ref{rem5} to obtain the following conclusion: if the $T$-periodic displacement-controlled loading $l(t)$ satisfies 
 (\ref{ex1_nonconstant}) and, for the $T$-periodic stress loading $h(t)$, one either has
\begin{eqnarray*}
& &   -\bar c_1^++\bar c^-_2-\bar c_3^+<a_1h_1(t)+a_2h_2(t)+a_3h_3(t)<\\
& & \hskip0.7cm<\min\left\{-\bar c_1^-+\bar c_2^--c_3^+,-\bar c_1^++\bar c_2^+-\bar c_3^+,-\bar c_1^++\bar c_2^--\bar c_3^-\right\},\quad t\in[0,T],
\end{eqnarray*}
or
\begin{eqnarray*}
&&\hskip-1.2cm  \max\left\{-\bar c_1^++\bar c_2^+-c_3^-,-\bar c_1^-+\bar c_2^--\bar c_3^-,-\bar c_1^-+\bar c_2^+-\bar c_3^+\right\}<\\
&&\hskip-0.5cm <a_1h_1(t)+a_2h_2(t)+a_3h_3(t)<
-\bar c_1^-+\bar c^+_2-\bar c_3^-,\quad t\in[0,T],
\end{eqnarray*}
then the stresses of springs of the elastoplastic system of Fig.~\ref{ex1fig} converge, as $t\to\infty,$  to a unique $T$-periodic regime that depends on $l(t)$ and $h(t)$, and doesn't depend on the initial state of the system.

\section{Conclusions}

\noindent We used Moreau sweeping process framework to analyze the asymptotic properties of quasistatic evolution of one-dimensional networks of elastoplastic springs (elastoplastic systems) under displacement-controlled and stress-controlled loadings. This type of elastoplastic systems covers, in particular, rheological models of materials science. We showed that displacement-controlled loading corresponds to parallel displacement of the moving polyhedron $C(t)$ of the respective sweeping process, but doesn't influence the shape of $C(t)$. 
We showed that it is the stress loading which is capable to change the shape of $C(t)$. Moreover, we proved that increasing the magnitude of the stress loading always makes $C(t)$ a simplex, 
if the number $q$ of displacement-controlled constraints is two less the number $n$ of nodes of the network ($q=n-2$).

\vskip0.2cm

\noindent The global asymptotic stability result established in this paper ensures convergence of the stresses of springs to a unique periodic solution (output) when the magnitude of the displacement-controlled loading is large enough and when the
normal vectors of any $d$ different facets of the moving polyhedron $C(t)$ are linearly independent.  Here $d$ is the dimension of the phase space of the  polyhedron $C(t)$, given by $d=m-n+q+1$, where $m$ is the number of springs, see (\ref{dimV}). The most natural example where such a property holds is when $C(t)$ is a simplex.
The paper, therefore,  puts the simplectic shape of $C(t)$ forward as a
{\it Discrete analogue of Drucker's postulate} (as far as materials science audience is concerned).

\vskip0.2cm \noindent Our theory can be viewed as an analogue of the high gain feedback stabilization of the classical control theory, see Isidori \cite[\S4.7]{isidori}. The high gain assumption of the control theory corresponds to our condition (\ref{bigforce}) on the magnitude of stress loading. Our assumption $q=n-2$ on the network of elastoplastic springs resembles the relative degree in control. 

\vskip0.2cm

\noindent The advantage of the proposed restriction  $\dim U=1$ is that it leads to simple analytic conditions (\ref{formulaprop3}) and (\ref{bigforce}) for the convergence of an elastoplastic system, which can be used for the design of elastoplastic systems that converge for the desired set of applied loadings.
Extending Theorem~\ref{mainthm} to the case where $\dim U>1$ is a doable task, but the respective inequality (\ref{bigforce}) transforms into a list of groups of inequalities, where the number of groups equals the number of selections of $\dim U$ from $m$ (equation (\ref{xi1}) gets replaced by the respective combinations of $\dim U$ equations). We don't see how such a condition can be useful in design of applied loadings, thus we stick to  $\dim U=1$.

\vskip0.2cm

\noindent The results of the paper can be extended to the case of dynamic evolution of elastoplastic systems with small inertia forces along the lines of Martins et al \cite{zamm}.

\vskip0.2cm

\noindent We like to think that the present paper opens a new room of  opportunities for researchers interested in applied analysis and control.
\begin{appendices}
\section{}\label{appen} {\bf Proof of Theorem~\ref{lemattract}} \textit{{ {\rm (}Massera-Krejci Theorem for sweeping processes with a moving set of the form {$C(t)=\cap_{i=1}^k (C_i+c_i(t))$}{\rm ).} }}

\vskip0.2cm
 
  \noindent We prove that 
every solution $x$ of sweeping process (\ref{swgen}), that is defined on $[0,\infty)$, satisfies
\begin{equation}\label{eq:periodicAttraction}
\lim_{t\to +\infty}\|x(t)-x^*(t)\|_0=0,
\end{equation}
where $x^*$ is a $T$-periodic solution of (\ref{swgen}).
\vskip0.2cm

\noindent Notice, that in case of $T-$periodic input the function $t\mapsto x(t+T)$ coincide with another solution of \eqref{swgen} originating from the point $x(T)$ at $t=0$. Due to monotonicity of $N_{C(t)}^0(x)$ in $x$ the distance 
$\|x(t+T)-x(t)\|_0$  is non-increasing(see e.g. \cite[Corollary~1]{kunze}) and
there exists
\begin{equation}
r=\lim_{t\to +\infty}\|x(t+T)-x(t)\|_0.
\label{eq:limitR}
\end{equation}


\noindent Since $x([0,\infty))$ is precompact, there is a subsequence $\{n_j\}_{j\in \mathbb{N}}\subset \mathbb{N}$ and a point $x^*_0$ such that
\begin{equation}\label{combinedwith}
\lim_{j\to \infty}\|x(n_j T)-x_0^*\|_0=0.
\end{equation}
Moreover, since each $x(n_j T)\in C(n_j T)=C(0)$ and $C(0)$ is closed we have $x^*_0\in C(0)$.
Let $x^*$ be a solution of (\ref{swgen}) with the initial condition $x^*(0)=x^*_0$. Consider  the functions 
$$x_j(t)=x(t+n_j T),\quad j\in\mathbb{N}.$$
Since $C(t)=C(n_j T+t)$, each function $x_j(t)$ is the solution of sweeping process (\ref{swgen}) with the initial condition $x_j(0)=x(n_j T)$.
The distance between solutions doesn't increase, so
for any $t\geqslant 0$,
\begin{equation}
 0\leqslant \|x_j(t)-x^*(t)\|_0\leqslant \|x(n_j T)-x^*_0\|_0,
\label{eq:x*isBounded}
\end{equation}
and using  (\ref{combinedwith}) we obtain
\eqref{eq:periodicAttraction}. Now it remains to prove that $x^*$ is $T$-periodic. Combining  (\ref{eq:x*isBounded}) and \eqref{eq:limitR} we get  
\begin{gather*}
r=\lim_{j\to\infty}\|x(t+n_j T+T)-x(t+n_j T)\|_0=
\|x^*(t+T)-x^*(t)\|_0.
\end{gather*}
Since $x^*$ and $t\mapsto x^*(t+T)$ are two solutions of sweeping process (\ref{swgen})  with the constant distance $r$ between them,  lemma~\ref{th:constDistance} yields 
\begin{equation}\label{onlynew}
\dot x^*(t)=\dot x^*(t+T),\quad t\ge 0.
\end{equation}
Thus, 
$$
x^*(\bar n T)-x^*_0=\int\limits_0^{\bar n T} \dot x^*(t) dt=\bar n\int\limits_0^\top \dot x^*(t) dt=\bar n(x^*(T)-x^*_0),\quad \bar n\in \mathbb{N},
$$
and so 
$\|x^*(\bar n T)-x^*_0\|=\bar nr,$ $\bar n\in \mathbb{N}.$ Since $t\mapsto x^*(t)$ is bounded, the latter is possible only when  $r=0$, i.e. when $x^*$ is $T$-periodic. The proof of the theorem is complete. \qed

\vskip0.2cm \noindent Similar to Theorem~\ref{lemattract} results are obtained in Henriquez \cite{Henriquez} (extension to Banach spaces) and in Kamenskii et al \cite{nahs} (extension to almost periodic solutions).
\section{Structure of the configuration space}
\label{appendB}
\begin{figure}[h]\center
\includegraphics{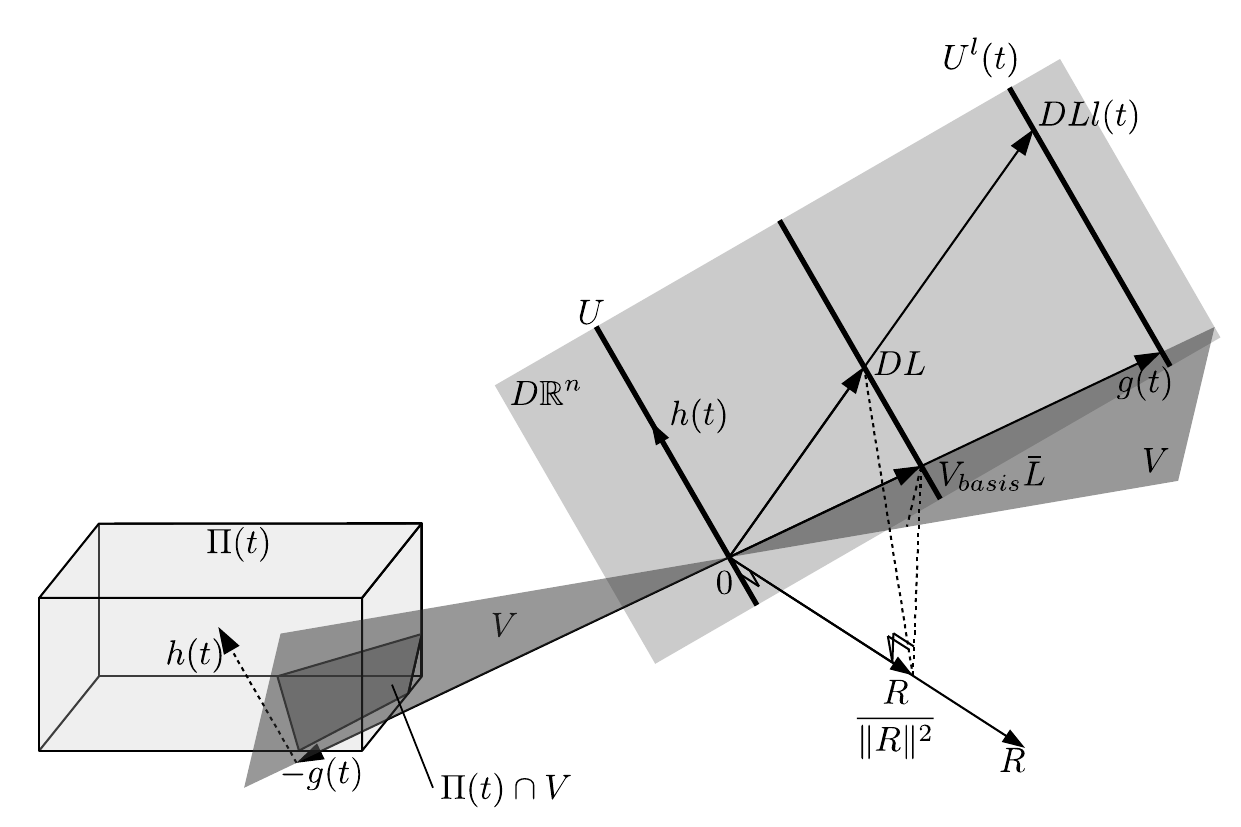}
\caption{\footnotesize The structure of the configuration space $\mathbb{R}^m$ when $q=1, {\rm dim}\, U=1, {\rm dim}\, V=2$. The right side of the figure shows how the vector $g(t)$ is obtained and the left side shows how the moving set $\Pi(t)\cap V$ is defined.} 
\label{fig:StructureOfSpace}
\end{figure}

To illustrate the structure of the configuration space $\mathbb{R}^m$ and construction of the variables used in Theorem~\ref{moreauthm}  we can plot a 3D  diagram (Figure \ref{fig:StructureOfSpace}) of a hypothetical situation when $m=3, q=1, {\rm rank}\, D =2$ whithout a  connection to any particular network of springs. In such a case we have ${\rm dim}\, U=1, {\rm dim}\, V=2$. Since $q=1$ the matrices $R$ and $L$ are single column-vectors and we illustrate the condition \eqref{barxi} on $L$ by showing that
\[
{\rm proj}_{R}(DL)=\frac{R}{\|R\|}\left\la\frac{R}{\|R\|}, DL\right\ra=\frac{R}{\|R\|^2} R^T DL=\frac{R}{\|R\|^2}
\]
\end{appendices}

\section*{Compliance with Ethical Standards}

\noindent {\bf Conflict of Interest:} The authors have no conflict of interest.

\bibliographystyle{plain}

\end{document}